\documentclass[11pt,reqno]{amsart}
\usepackage{euscript,amssymb,amsbsy,mathabx}
\usepackage[arrow,matrix,curve]{xy}
\usepackage{array,longtable,graphicx,enumitem,xcolor}
\usepackage{a4wide}
\definecolor{darkblue}{RGB}{40,40,85}
\usepackage[colorlinks=true,linkcolor=black,urlcolor=darkblue,citecolor=darkblue]{hyperref}

\numberwithin{equation}{section}
\numberwithin{figure}{section}

\newcounter{msct}[section]\renewcommand{\themsct}{\thesection.\arabic{msct}}
\newenvironment{m-theorem}{\vskip3pt\refstepcounter{msct}\trivlist \itemindent 0pt %
\item[\hskip\labelsep\bf Theorem \themsct]\it\ignorespaces}{\endtrivlist\vskip2pt}
\newenvironment{m-proposition}{\vskip3pt\refstepcounter{msct}\trivlist \itemindent 0pt %
\item[\hskip\labelsep\bf Proposition \themsct]\it\ignorespaces}{\endtrivlist\vskip2pt}
\newenvironment{m-corollary}{\vskip3pt\refstepcounter{msct}\trivlist \itemindent 0pt %
\item[\hskip\labelsep\bf Corollary \themsct]\it\ignorespaces}{\endtrivlist\vskip2pt}
\newenvironment{m-lemma}{\vskip3pt\refstepcounter{msct}\trivlist \itemindent 0pt %
\item[\hskip\labelsep\bf Lemma \themsct]\it\ignorespaces}{\endtrivlist\vskip2pt}
\newenvironment{m-definition}{\vskip3pt\refstepcounter{msct}\trivlist \itemindent 0pt %
\item[\hskip\labelsep\bf Definition \themsct]\ignorespaces}{\endtrivlist\vskip2pt}
\newenvironment{m-notation}{\vskip3pt\refstepcounter{msct}\trivlist \itemindent 0pt %
\item[\hskip\labelsep\bf Notation \themsct]\ignorespaces}{\endtrivlist\vskip3pt}
\newenvironment{m-example}{\vskip3pt\refstepcounter{msct}\trivlist \itemindent 0pt %
\item[\hskip\labelsep\bf Example \themsct]\ignorespaces}{\endtrivlist\vskip3pt}
\newenvironment{m-remark}{\vskip3pt\refstepcounter{msct}\trivlist \itemindent 0pt %
\item[\hskip\labelsep\bf Remark \themsct]\ignorespaces}{\endtrivlist\vskip3pt}
\newenvironment{m-question}{\vskip3pt\refstepcounter{msct}\trivlist \itemindent 0pt %
\item[\hskip\labelsep\bf Question.]\ignorespaces}{\endtrivlist\vskip3pt}
\newenvironment{thm-nono}[1]{\vskip3pt\trivlist \itemindent 0pt %
\item[\hskip\labelsep\bf Theorem #1]\it\ignorespaces}{\endtrivlist\vskip3pt}
\newenvironment{conj-nono}[1]{\vskip3pt\trivlist \itemindent 0pt %
\item[\hskip\labelsep\bf Conjecture #1]\it\ignorespaces}{\endtrivlist\vskip3pt}
\newenvironment{m-thank}{\vskip3pt\trivlist \itemindent 0pt %
\item[\hskip\labelsep\it Acknowledgments]\ignorespaces}{\endtrivlist\vskip3pt}
\newenvironment{m-proof}{\vskip2pt\trivlist \itemindent 0pt %
\item[\hskip\labelsep\it Proof.]\ignorespaces}{\hfill$\Box$\endtrivlist\vskip3pt}%
\newenvironment{m-asmp}{\vskip3pt\trivlist \itemindent 0pt %
\item[\hskip\labelsep\bf Assumption.]\ignorespaces}{\hfill\endtrivlist\vskip3pt}%

\newcommand{\bibauth}[2]{\textrm{{#1}~{#2}}}
\newcommand{\bibtitl}[1]{\textit{#1}}
\newcommand{\bibjnyp}[4]{\textrm{#1} \textbf{#2} (#3), #4}
\newcommand{\bibinbook}[3]{In: \textrm{#1}\textrm{, #2}\textrm{, #3}}
\newcommand{\bibbook}[4]{\textit{#1}. {#2} {#3}, {#4}}

\let\lar\longrightarrow

\let\hra\hookrightarrow
\let\mt\mapsto

\font\tenmsa=msam10 %
\newcommand\hdashpiece{%
{\vrule height2.75pt depth-2.35pt width2.3pt \kern1.7pt}}%
\newcommand\hdashpieces{%
{\hdashpiece\hdashpiece\hdashpiece\hdashpiece}}%
\let\dashto\dashrightarrow
\newcommand\dashar{\mathrel{%
\hdashpieces\kern-0.4pt\hbox{\tenmsa K}}}%

\let\euf\EuScript 
\let\cal\mathcal
\let\mbb\mathbb

\let\mfrak\mathfrak
\let\bsymb\boldsymbol 
\DeclareFontFamily{OT1}{rsfs}{}
\DeclareFontShape{OT1}{rsfs}{n}{it}{<->rsfs10}{}
\DeclareMathAlphabet{\crl}{OT1}{rsfs}{n}{it}

\newcommand\uset[2]{{\disp\mathop{\mbox{$#2$}}_{#1}}}
\newcommand\oset[2]{{\disp\mathop{\mbox{$#2$}}^{#1}}}
\newcommand\ouset[3]{{\oset{#2}{\uset{#1}{#3}}}}
\let\ovl\overline

\let\unbar\underbar

\let\tld\tilde
\let\wtld\widetilde

\let\nit\noindent

\let\disp\displaystyle
\let\srel\stackrel
\let\vphi\varphi
\let\eps\epsilon
\let\veps\varepsilon

\newcommand\lran[1]{{\langle #1\rangle}}

\newcommand\ort{\mathrel{{\vrule width4.0pt height0.4pt depth0pt
                \vrule width0.4pt height6.0pt depth0pt\,}}}
\newcommand\rd{{\rm d}} 

\newcommand\Aut{\operatorname{\textrm{Aut}\kern1pt}}
\newcommand\cAut{\operatorname{\mathcal{A}\kern-1pt\textit{ut}\kern1pt}}
\newcommand\End{\operatorname{\rm{End}\kern1pt}}
\newcommand\cEnd{\operatorname{\mathcal{E}\kern-1pt\textit{nd}\kern1pt}}

\newcommand\Hom{\mathop{\rm Hom}\nolimits}
\newcommand\cHom{\operatorname{\mathcal{H}\kern-1pt\textit{om}\kern1pt}}
\newcommand\Img{{\rm Im}}
\newcommand\Ker{{\rm Ker}}

\newcommand\NS{\mathop{\rm NS}\nolimits}
\newcommand\Pic{\mathop{\rm Pic}\nolimits}
\newcommand\Proj{\mathop{\rm Proj}\nolimits}

\newcommand\Spec{\mathop{\rm Spec}\nolimits}

\newcommand\Sym{\mathop{\rm Sym}\nolimits}
\newcommand\invq{{\slash\kern-2.5pt\slash}}

\newcommand\rk{{\rm rk}}

\newcommand\Grs{{\rm Gr}}
\newcommand\spGrs{{\rm sp{\hdashpiece\!}Gr}}
\newcommand\oGrs{{\rm o{\hdashpiece\!}Gr}}

\let\lda\lambda

\let\si\sigma
\let\Si\Sigma

\let\sm\setminus

\newcommand\bbC{{\mbb C}}

\newcommand\bbk{\mbox{\rm I\kern-1.5pt k}}
\newcommand\sbbk{\hbox{\scriptsize I{\kern-.8pt}k}}

\newcommand\bone{{1\kern-0.57ex\rm l}}

\newcommand\eA{{\euf A}}

\newcommand\cE{{\cal E}}
\newcommand\eE{{\euf E}}
\newcommand\cF{{\cal F}}
\newcommand\eF{{\euf F}}
\newcommand\eG{{\euf G}}
\newcommand\eI{{\euf I}}
\newcommand\eJ{{\euf J}}

\newcommand\cL{{\cal L}}
\newcommand\eL{{\euf L}}

\newcommand{\eM}{{\euf M}}

\newcommand\eN{{\euf N}}
\newcommand\cO{{\cal O}}
\newcommand\eO{{\euf O}}

\newcommand\cY{{\cal Y}}

\newcommand\eT{{\euf T}}

\newcommand\codim{{\rm codim}}

\newcommand\Ad{{\rm Ad}}
\newcommand\Bl{{\rm Bl}}

\newcommand\SO{{\rm SO}}

\newcommand\Sp{{\rm Sp}}

\renewcommand\det{{\rm det}}

\let\ges\geqslant
\let\les\leqslant

\newcommand\Supp{\mathop{\rm Supp}\nolimits}

\newcommand\res{\mathop{\rm res}\nolimits}
\let\surj\twoheadrightarrow

\newcommand\pos{{>0}}
\newcommand\apos{{\,\gtrsim0}}

\newcommand\diag{\mathop{\rm diag}\nolimits}

\newcommand{\BB}{{\rm BB}}

\newcommand{\cd}{\mathop{\rm cd}\nolimits}

\newcommand{\bnom}[2]{\genfrac{(}{)}{0pt}{}{#1}{#2}}
\newcommand{\cst}[1]{\mathop{\rm ct}^{#1}\nolimits}
\newcommand{\cX}{{\mathcal X}}
\newcommand{\Hilb}{\mathop{\rm Hilb}\nolimits}
\newcommand{\kk}{{\Bbbk}}
\newcommand{\reg}{\mathop{\rm reg}\nolimits}
\newcommand{\trdeg}{\mathop{\rm trdeg}\nolimits}
\newcommand{\rst}{{\upharpoonright}}
\newcommand{\bz}{{\textbf{z}}}

\author{Mihai Halic}
\address{}
\email{mihai.halic@gmail.com}


\begin{document}

\title[$q$-ample subvarieties]{Subvarieties with $q$-ample normal bundle\\ and $q$-ample subvarieties}

\begin{abstract}
The goal of this article is twofold. On one hand, we study the subvarieties of projective varieties which possess partially ample normal bundle; we prove that they are G2 in the ambient space. This generalizes results of Hartshorne and B\u{a}descu-Schneider. We work with the cohomological partial ampleness introduced by Totaro. 

On the other hand, we define the concept of a partially ample subvariety, which generalizes the notion of an ample subvariety introduced by Ottem. We prove that partially ample subvarieties enjoy the stronger G3 property. Moreover, we present an application to a connectedness problem posed by Fulton-Hansen and Hartshorne. 

The results are illustrated with examples.
\end{abstract}

\maketitle

\section*{Introduction}

While searching for an appropriate concept of amplitude for higher codimensional subvarieties, Hartshorne investigated the geometric and cohomological properties of pairs $(X,Y)$ consisting of a projective scheme $X$ and a local complete intersection subscheme $Y$ with ample normal bundle. On one hand, $Y\subset X$ is G2 that is, the formal completion $\hat{X}_Y$ determines an \'etale neighbourhood of $Y$; furthermore, subvarieties of projective spaces are actually G3 (cf. \cite{hart-cdav,hart-as,hir+mats}). On the other hand, the cohomology groups of coherent sheaves on the complement $X\sm Y$ are finite dimensional and vanish, above appropriate degrees. 

It turns out that the assumption about the ampleness of the normal bundle can be weakened. It suffices to require either a Hermitian metric with partially positive curvature (cf. \cite{griff,comm+graut}) or partial ampleness in the sense of Sommese  (cf. \cite{hirw,bad+schn}). A comprehensive reference to the algebraic approach to the problem is B\u{a}descu's book \cite{bad}. 

Ample subvarieties of projective varieties were defined by Ottem in \cite{ottm}, based on Totaro's work \cite{tot} on cohomological ampleness. Their normal bundle is ample, but not conversely. Ample subvarieties enjoy pleasant features: they are G3 in the ambient space and the cohomological dimension of their complement is the smallest possible (the co-dimension minus one). 

In a different framework, the author introduced the weaker notion of a $q$-ample subvariety, generalizing the ample subvarieties defined by Ottem (cf. \cite{hlc-qvb}). One shows that the $q$-ampleness of a subvariety $Y\subset X$ consists of two conditions:
\begin{itemize}[leftmargin=5ex]
\item[--] 
a global one, which is an upper bound on the cohomological dimension of the complement $X\sm Y$;
\item[--] 
a local one, the $q$-ampleness of the normal bundle.
\end{itemize}

In general, it is difficult to estimate the cohomological dimension of the complement of a subvariety. The issue has been intensively investigated (cf. \cite{hart-as,ogus,lyub-licht,falt-homog}). 

Concerning the second condition, Hartshorne initiated the study of subvarieties with ample normal bundle (cf. \cite{hart-cdav,hart-as}). In \cite{bad+schn}, B\u{a}descu-Schneider proved the G2 property for subvarieties with \emph{globally generated}, partially ample normal bundle (in the sense of Sommese \cite{soms}). The result implies that generating subvarieties of rational homogeneous and abelian varieties are G3. Their approach essentially uses the global generation of the normal bundle, in order to reduce the problem to \cite{hart-cdav}. Faltings proved in \cite{falt-homog} that the G3 property holds for low codimensional subvarieties of rational homogeneous varieties and also estimated the cohomological dimension of their complement.

To our knowledge, the properties of subvarieties with $q$-ample normal bundle have not been investigated yet. We emphasize that we refer to the cohomological partial ampleness introduced by Arapura-Totaro \cite{arap,tot}, which is less restrictive than the partial ampleness of Sommese. This is the main motivation for our study. 

The advantage of the cohomological ampleness is that of being a numerical condition. We show that the G2 property holds for subvarieties whose co-normal bundle is not pseudo-effective. This yields numerous \emph{new examples} of subvarieties with partially ample, but neither ample nor globally generated, normal bundle: \textit{e.g.} sufficiently generic subvarieties of uniruled varieties. The ubiquity of these objects is, in our opinion, a strong motivation to systematically study their properties. 
We highlight the main results obtained in this article. 
\begin{thm-nono}{} 
Let $X$ be a non-singular irreducible projective variety, defined over an algebraically closed field of characteristic zero, and let $Y$ be a local complete intersection subvariety. Then the following statements hold: 
\begin{enumerate}[leftmargin=5ex]
\item[\rm(i)] (cf.~\ref{thm:finite-dim}, \ref{cor:G2}) 
If $Y$ is connected and its normal bundle in $X$ is $(\dim Y-1)$-ample, then $Y$ is G2 and the cohomology groups of coherent sheaves on $X\sm Y$ are finite dimensional, above a certain degree.
\item[\rm(ii)] (cf.~\ref{thm:G3}) 
If $Y$ is $(\dim Y-1)$-ample, then it is G3. In particular, it holds: 
\begin{enumerate}
\item[\rm(a)] (cf.~\ref{thm:badescu2}) 
Sufficiently general, movable subvarieties of rationally connected varieties are G3. 
\item[\rm(b)] (cf.~\ref{cor:DX}) 
Suppose $X$ is almost homogeneous for the action of a linear algebraic group $G$, with open orbit $O$. If $\codim_X(X\sm O)\ges 2$ and the stabilizer of a point in $O$ contains a Cartan subgroup of $G$, then the diagonal of $X\times X$ is G3.
\end{enumerate}
\end{enumerate}
\end{thm-nono}
Moreover, we prove the G3 property for movable subvarieties of minimal Mori dream spaces (cf.~\ref{cor:mori-g3}) and of the strongly (very) movable subvarieties introduced by Voisin (cf.~\ref{thm:strg-mov}). Furthermore, we obtain several results concerning the behaviour under natural operations, \textit{e.g.} intersections, products, pre-images. 

We discuss a conjecture of Fulton-Hansen \cite{fult+hans} and Hartshorne \cite{hart-as}, concerning the connectedness of the intersection of two subvarieties (of an ambient space $X$) with ample normal bundle. The conjecture and various versions hold for $X$ a projective space, a flag variety, a homogeneous space, or products of such (cf. \cite{fult+hans,hans,falt-homog,debr2,bad-debarre}). 
The approach inhere applies to arbitrary ambient varieties. We impose partial ampleness to one subvariety---not only to its normal bundle---and we analyse how is this property preserved by intersections. Our contribution to this topic is the following (see \textit{loc. cit.} for details): 

\begin{thm-nono}{{\rm(cf. \ref{thm:fh})}} 
Let $X$ be a projective variety, $Y$ a lci subvariety, and $V\srel{f}{\to} X$ a morphism, with $V$ projective, irreducible. Suppose $q:=\dim f(V)-\codim_XY>0$ and $Y\subset X$ is $q$-ample. Then the following statements hold: 
\begin{enumerate}[leftmargin=5ex]
\item[\rm(i)]
If the Stein factorization $\bar V=\Spec(f_*\eO_V)$ of $f$ is a Cohen-Macaulay variety, then $f^{-1}(Y)$ is non-empty and connected.
\item[\rm(ii)] 
Assume that $f(V)$ is smooth, $Y\cap f(V)$ is lci in $f(V)$. Then $f^{-1}(Y)\subset V$ is G3.
\end{enumerate}
\end{thm-nono}

A substantial part is dedicated to illustrate the general theory with examples. It is pleasing that partially ample subvarieties naturally occur in various contexts: 
\begin{enumerate}[leftmargin=5ex]
\item[(i)] 
almost homogeneous, more generally, rationally connected varieties (cf. {\S}\ref{ssct:subvar-homog}); 
\item[(ii)] 
zero loci of sections in globally generated vector bundles (cf. {\S}\ref{ssct:glob-gen}); 
\item[(iii)]
Bialynicki-Birula decompositions, corresponding to $G_m$-actions (cf. {\S}\ref{ssct:fixed}). 
\end{enumerate}

The article is divided in two. In the first part, we follow \cite{hart-cdav} to study the properties of subvarieties with partially ample normal bundle; their main property is that of being G2 in the ambient space. We apply the theory to the case of uniruled varieties and generalize existing results which typically hold for complete intersections or for subvarieties of homogeneous spaces (cf. \ref{thm:mov}, \ref{cor:diag-G2}). 

In the second part, we define and study partially ample subvarieties, in a similar vein to \cite{ottm}. One of their essential features is that of being G3; the proof of this fact depends on the G2-result obtained in the first part. We adapt some of the results in \cite{bad} to deduce the G3 property for subvarieties of almost homogeneous and of rationally connected varieties. The article ends with explicit computations in several cases: zero loci of sections in vector bundles and sources of actions of the multiplicative group.

\setcounter{tocdepth}{1}\tableofcontents
\section{Preliminary notions}\label{sct:def}

\begin{m-notation}\label{XYN} 
We work over an algebraically closed field $\kk$ of characteristic zero. Throughout the article, $\mfrak X$ stands for a connected, noetherian formal scheme, regular and projective over $\kk$; $X$ stands for an irreducible projective variety over $\kk$. 

We denote by $Y$ either a subscheme of definition of $\mfrak X$---by assumption, it is projective---, or a closed subscheme of $X$. Let $\dim Y$ be the maximal dimension of its irreducible components; if $Y\subset X$, $\codim(Y):=\dim X-\dim Y$. We assume that all the components of $Y$ are at least $1$-dimensional, so $\dim Y\ges1$.

Furthermore, let $\eI_Y\subset\cal O_{\mfrak X}$ (resp. $\subset\eO_X$) be the sheaf of ideals defining $Y$. For $a\ges0$, let $Y_{a}$ be the subscheme defined by the $\eI_{Y}^{a+1}$. 
The formal completion of $X$ along $Y$ is defined as $\hat X_Y:=\disp\varinjlim Y_a$. If $X$ is non-singular in a neighbourhood of $Y$, then $\hat X_Y$ is regular and projective. Details can be found in \cite[Section II.9]{hart-ag}. 

In the case where $Y$ is a locally complete intersection---\emph{lci} for short---in $\mfrak X$, we denote its normal sheaf by $\eN=\eN_Y:=(\eI_Y/\eI_Y^2)^\vee$; it is locally free of rank $\nu$ on $Y$. The structure sheaves of the various thickenings $Y_a$ fit into the exact sequences:
\begin{equation}\label{eq:Ya}
0\to\Sym^a(\eN^\vee)\to\eO_{Y_a}\to\eO_{Y_{a-1}}\to0,\;\forall a\ges1.
\end{equation}

For a coherent sheaf $\eG$, we denote $h^t(\eG):=\dim_{\kk}H^t(\eG)$; for a field extension $K\hra K'$, $\trdeg_KK'$ is the transcendence degree; $\cst{A,B,...}$ stands for a real constant depending on the quantities $A,B,\dots$. A \emph{line (resp. vector) bundle} is an \emph{invertible (resp. locally free) sheaf}.

Suppose $Y$ is connected; let $K(\hat X_Y)$ be the field of formal rational functions on $X$ along $Y$. We recall the following terminology due to Hironaka-Matsumura (cf. \cite{hir+mats}):
\begin{itemize}[leftmargin=5ex]
\item[--] $Y$ is G1 in $X$, if $K(\hat X_Y)=\kk$;
\item[--] $Y$ is G2 in $X$, if $K(X)\hra K(\hat X_Y)$ is finite;
\item[--] $Y$ is G3 in $X$, if $K(X)\hra K(\hat X_Y)$ is an isomorphism.
\end{itemize}
\end{m-notation}

\subsection{Cohomological $q$-ampleness}\label{ssct:cohom-q} 
This notion was introduced by Arapura and Totaro. 

\begin{m-definition}\label{def:q-line}
Let $Y$ be a projective scheme. 

\nit{\rm(i)} (cf. \cite[Theorem 7.1]{tot}) 
An invertible sheaf $\eL$ on $Y$ is called \emph{$q$-ample} if, for any coherent sheaf $\eG$ on $X$, holds: 
\begin{equation}\label{eq:qL}
\exists\,\cst{\eG}\;\forall\,a\ges \cst{\eG}\;\forall\,t>q,\quad H^t(Y,\eG\otimes\eL^{a})=0.
\end{equation}
It's enough to verify this property for $\eG=\eA^{-k}, k\ges1$, where $\eA\in\Pic(Y)$ is a fixed, but otherwise arbitrary, ample line bundle (cf. \cite[Theorem 6.3]{tot}).

\nit{\rm(ii)} (cf. \cite[Lemma 2.1, 2.3]{arap}) 
A locally free sheaf $\eE$ on $Y$ is \emph{$q$-ample} if $\eO_{\mbb P(\eE^\vee)}(1)$ on $\mbb P(\eE^\vee):=\Proj(\Sym^\bullet_{\eO_Y}\eE)$ is $q$-ample. This is equivalent saying that, for any coherent sheaf $\eG$ on $Y$, there is $\cst{\eG}>0$ such that 
\begin{equation}\label{eq:qG}
H^t(Y,\eG\otimes\Sym^a(\eE))=0,\;\forall t>q,\; \forall a\ges \cst{\eG}.
\end{equation}
The \emph{$q$-amplitude of $\eE$}, denoted $q^\eE$, is the smallest integer $q$ which satisfies this property. 
\\ We remark that $\eE$ is $q$-ample if and only if $\eE_{Y_{\rm red}}$ is $q$-ample (cf. \cite[Corollary 7.2]{tot}). 
\end{m-definition}

The $q$-amplitude enjoys \emph{uniformity} and \emph{sub-additivity} properties, which we recall.

\begin{m-theorem}\label{thm:unif-q}\quad 
{\rm(i) (cf. \cite[Theorem 6.4, 7.1]{tot})}
Let $Y$ be a projective scheme, $\eA,\eL\in\Pic(Y)$. We assume that $\eA$ is sufficiently ample---Koszul-ample, cf. \cite[pp. 733]{tot}---, and $\eL$ is $q$-ample. Then there are constants $\cst{\eA,\eL}_1,\cst{\eA,\eL}_2>0$, such that for any coherent sheaf $\eG$ on $Y$ holds:
$$
H^t(Y,\eG\otimes\eL^a)=0,\;\;\forall t>q,\,
\forall a\ges\cst{\eA,\eL}_1\cdot\reg^\eA(\eG)+\cst{\eA,\eL}_2.
$$
Here $\reg^\eA(\eG)$ stands for the regularity of $\eG$ with respect to $\eA$. 

\nit{\rm(ii) (cf. \cite[Theorem 3.4]{tot})} 
If $H^0(\eO_Y)=\kk$ then, for any locally free sheaf $\eE$ and coherent sheaf $\eG$ on $Y$, holds: $\;\reg^\eA(\eE\otimes\eG)\les\reg^\eA(\eE)+\reg^\eA(\eG)$. 
\end{m-theorem}

\begin{m-theorem}\label{thm:subadd-q}{\rm(cf. \cite[Theorem 3.1]{arap})}
Let $0\to\eE_1\to\eE\to\eE_2\to 0$ be an exact sequence of locally free sheaves on the projective scheme $Y$. Then it holds: $\;q^\eE\les q^{\eE_1}+q^{\eE_2}$. 
\end{m-theorem}

For products one has a better estimate. 

\begin{m-lemma}\label{lm:q-prod}
Let $X_1,X_2$ be irreducible projective varieties and $\eE_1,\eE_2$ be locally free sheaves on them which are $q^{\eE_1}$-, $q^{\eE_2}$-ample, respectively. Let $\eE_1\boxplus\eE_2$ be the direct sum of their pull-backs to $X_1\times X_2$. Then we have: 
$\;q^{\eE_1\boxplus\eE_2}\les\max\{q^{\eE_1}+\dim X_2,q^{\eE_2}+\dim X_1\}.$ 
\end{m-lemma}

\begin{m-proof}
It suffices to take $\eG$ in \eqref{eq:qG} of the form $(\eA_1\boxtimes\eA_2)^{-1}$, with $\eA_1,\eA_2$ ample line bundles on $X_1,X_2$, respectively. Then holds:
$$
\begin{array}{l}
H^t\big(X_1\times X_2,(\eA_1^{-1}\boxtimes\eA_2^{-1})\otimes\Sym^a(\eE_1\boxplus\eE_2)\big)
\\[1.5ex]
=
\uset{\genfrac{}{}{0pt}{}{t_1+t_2=t,}{a_1+a_2=a}}{\bigoplus}
H^{t_1}\big(X_1,\eA_1^{-1}\otimes\Sym^{a_1}(\eE_1)\big)\otimes
H^{t_2}\big(X_2,\eA_2^{-1}\otimes\Sym^{a_2}(\eE_2)\big). 
\end{array}
$$
It remains to apply the definition, for $t\ges\max\{q^{\eE_1}+\dim X_2,q^{\eE_2}+\dim X_1\}$. 
\end{m-proof}

\subsection{$(\dim Y-1)$-ample vector bundles on $Y$}\label{ssct:y-1}

Subvarieties $Y\subset X$ with $(\dim Y-1)$-ample normal bundle play an essential role in this article. Here we give a numerical characterization of this property, analogous to Totaro's result for invertible sheaves. 

\begin{m-proposition}\label{prop:y-1}{\rm(cf. \cite[Theorem 9.1]{tot})}
Let $\eE$ be a locally free sheaf on an irreducible projective variety $Y$. The following statements are equivalent: 
\begin{enumerate}[leftmargin=5ex]
\item[\rm(i)] 
$\eE$ is $(\dim Y-1)$-ample. 

\item[\rm(ii)] 
$\eO_{\mbb P(\eE)}(1)$ is not pseudo-effective, where $\mbb P(\eE):=\Proj(\Sym^\bullet\eE^\vee)$;\\ 
in this case, we say that \emph{$\eE^\vee$ is not pseudo-effective}. 

\item[\rm(iii)]  
There is a dominant morphism $\vphi:C_S\to Y$, with $S$ affine and $C_S$ an integral curve over $S$, such that the following conditions are satisfied: 
\begin{enumerate}
\item[\rm(1)] 
$\vphi^*\eE$ admits a line sub-bundle $\eM$ which is relatively ample for $C_S\to S$; 
\item[\rm(2)] 
Denote by $C_{S,y}$ the curves passing through the general point $y\in Y$. Then the corresponding lines $\eM_{y}\subset\eE_y$ are movable in $\mbb P(\eE_y)$. 
\end{enumerate}
\end{enumerate}
\end{m-proposition}

\begin{m-proof}
Let $\eO_Y(1)$ be an ample line bundle on $Y$. Denote by $\omega_Y$ the dualizing sheaf of $Y$ (cf. \cite[Prop. III.7.5]{hart-ag}); it is torsion free of rank one and there is $c>0$ such that:  
$$
\omega_Y\subset\eO_Y(c), \;\eO_Y(-c)\subset\omega_Y\quad\text{(cf. \cite[\S9]{tot})}.
$$
The $(\dim Y-1)$-ampleness means:  
$H^0(Y,\omega_Y\otimes\eL\otimes\Sym^a\eE^\vee)=0,$ $\forall\eL\in\Pic(Y),$ $a>\cst{\eL}$, which is equivalent to $H^0(Y,\eM\otimes\Sym^a\eE^\vee)=0,$ $\forall\eM\in\Pic(Y),\, a>\cst{\eM}$, and to: 
$$
H^0(\mbb P(\eE),\eM\otimes\eO_{\mbb P(\eE)}(a))=0,\,\forall\eM\in\Pic(\mbb P(\eE)),\forall a>\cst{\eM}.
$$
But the last condition is the $(\dim\mbb P(\eE)-1)$-ampleness of $\eO_{\mbb P(\eE)}(-1)$ and (i)$\Leftrightarrow$(ii) follows. 

The equivalence (ii)$\Leftrightarrow$(iii) is a direct consequence of the duality theorem between the pseudo-effective cone and the movable cone \cite[Theorem 0.2]{bdpp}. 
\end{m-proof}

Note that the notion of a \emph{pseudo-effective vector bundle} used in \cite[\S7]{bdpp} is more restrictive. It also requires that the projection of the non-nef locus of $\eO_{\mbb P(\eE)}(1)$ does not cover $Y$.


\subsection{$q$-positivity}

\begin{m-proposition}{\rm (cf. \cite[Proposition 1.7]{soms}).}\label{prop:q2}
For a \emph{globally generated}, locally free sheaf $\eE$ on $Y$, the following statements are equivalent: 
\begin{enumerate}[leftmargin=5ex]
\item[\rm(i)] 
$\eE$ is $q$-ample (cf. Definition \ref{def:q-line}); 
\item[\rm(ii)] 
The fibres of the morphism $\mbb P(\eE^\vee)\to |\eO_{\mbb P(\eE^\vee)}(1)|$ are at most $q$-dimensional. 
\end{enumerate}
We say that $\eE$ is \emph{Sommese-$q$-ample} if it satisfies any of these conditions.
\end{m-proposition}

\begin{m-definition}\label{def:q-pos}(cf. \cite{andr+grau,dps}) 
Suppose $X$ is a smooth, complex projective variety. A line bundle $\eL$ on $X$ is \emph{$q$-positive}, if it admits a Hermitian metric whose curvature is positive definite on a subspace of $\eT_{X,x}$ of dimension at least $\dim X-q$, for all $x\in X$;  equivalently, the curvature has at each point $x\in X$ at most $q$ negative or zero eigenvalues. 
\end{m-definition}

\begin{m-theorem}\label{thm:mats-bott}\quad 
\nit{\rm(i)}{\rm(cf. \cite[Theorem 1.4]{mats})} 
Assume $\eE$ is globally generated. Then it holds: 
$$\eE\text{ is Sommese-$q$-ample}\;\Leftrightarrow\;
\eO_{\mbb P(\eE^\vee)}(1)\text{ is $q$-positive.}$$ 

\nit{\rm(ii)}{\rm(cf. \cite{bott-lefz,ottm})} 
Let $\eL\in\Pic(X)$ be $q$-positive and $Y\in|\eL|$ a smooth divisor. Then it holds: 
$$
H^t(X;\mbb Z)\to H^t(Y;\mbb Z)\text{ is } 
\bigg\{\begin{array}{rl}
\text{an isomorphism, for}&t\les \dim X-q-2;
\\[1.25ex] 
\text{surjective, for}&t=\dim X-q-1.
\end{array}
$$
\end{m-theorem}


\section*{\textbf{Part~I:\hspace{2ex}Subvarieties with \textit{q}-ample normal bundle}}

Hartshorne investigated in \cite{hart-cdav} the cohomological properties of subvarieties with ample normal bundle and of their complements. Here we generalize a number of his results---those in \textit{op.cit.}, Sections 5,~6---to subvarieties with cohomologically $q$-ample normal bundle. The difficulty to overcome is that several statements, especially those needed to deduce the G2 property, are proved for \emph{curves} with ample normal bundle; the general case is obtained by induction. For this reason, we must reprove the statements involving curves. 

B\u{a}descu and Schneider \cite{bad+schn} generalized Hartshorne's results to subvarieties whose normal bundle is Sommese-$q$-ample. They assume that the normal bundle is globally generated, because the proof is based on dimensional reduction. Hence their applications concern mainly subvarieties of homogeneous spaces and abelian varieties.


\section{Finite dimensionality results}\label{sct:finite-dim}

In this section we follow \cite[Section 5]{hart-cdav}. 

\begin{m-theorem}\label{thm:finite-dim}
{\rm(cf. \cite[Theorem 5.1, Corollary 5.4]{hart-cdav})} 
Let the situation be as in \ref{XYN}. Assume that $Y$ is lci and $\eN$ is $q^\eN$-ample. For any locally free sheaf $\cal{F}$ on $\mfrak{X}$ (of finite rank), holds: 
$$H^{t}(\mathfrak{X},\cal{F})\text{ is finite dimensional, for }t<\dim Y-q^\eN.$$ 

\nit In particular, if $q^\eN\leqslant\dim Y-1$ and $\mfrak{X}$ is connected, then $H^{0}(\mfrak{X},\cO_{\mfrak{X}})=\kk$.
\end{m-theorem}

\begin{m-proof}
Since $\disp H^t(\mfrak X,\cal F)=\varprojlim H^t(Y_a,\cal F\otimes\eO_{Y_a})$, it is enough to prove that the sequence eventually becomes stationary. For $\eF:=\cal F\otimes\eO_Y$, the exact sequences \eqref{eq:Ya} yield:
$$
\dots\to H^t(Y,\eF\otimes\Sym^a(\eN^\vee))
\to H^t(Y_a,\cal F\otimes\eO_{Y_a})
\to H^t(Y_{a-1},\cal F\otimes\eO_{Y_{a-1}})
\to\dots
$$
The lci condition implies that $Y$ is Gorenstein. The Serre duality and the $q^\eN$-ampleness of $\eN$ imply the vanishing of the leftmost term, for $a\gg0$. For the second statement, observe that $H^{0}(\mfrak{X},\cO_{\mfrak{X}})$ is an integral domain and also a finite dimensional $\kk$-algebra.
\end{m-proof}

\begin{m-corollary}\label{cor:finite-dim}
{\rm(cf. \cite[Corollary 5.3, 5.5]{hart-cdav})}  
\begin{enumerate}[leftmargin=5ex]
\item[\rm(i)] 
Let the situation be as in \ref{thm:finite-dim}. Suppose $\cL$ is an invertible sheaf on $\mfrak{X}$, such that its restriction to $Y$ is $q^\cL$-ample. Then holds: 
$$
H^{t}(\mfrak{X},\cal{F}\otimes\cL^{-b})=0,
\text{ for }t<\dim Y-(q^\eN+q^\cL),\,b\gg0.
$$
\item[\rm(ii)] 
Let $X$ be a non-singular projective scheme over $\kk$ and $Y$ a closed lci subscheme whose normal bundle is $q^\eN$-ample; let $\eL$ be a $q^\eL$-ample, invertible sheaf on $X$. Then the following statements hold, for all coherent sheaves $\eG$ on $X\setminus Y$: 
$$
\begin{array}{ccl}
(1)&H^{t}(X\setminus Y,\eG)\text{ is finite dimensional, }
&
t\ges\dim X-\dim Y+q^\eN,
\\[1ex] 
(2)&H^{t}(X\setminus Y,\eG\otimes\eL^b)=0,
&
t\ges\dim X-\dim Y+q^\eN+q^\eL,\;b\gg0.
\end{array}
$$
\end{enumerate}
\end{m-corollary}

\begin{m-proof}
(i) Denote $\eF:=\cal{F}\otimes\eO_{Y}$, $\eL:=\cL\otimes\eO_{Y}$; by \eqref{eq:Ya}, it is enough to show:
$$\;
H^{t}(Y,\omega_Y\otimes\eF^\vee\otimes\Sym^{a}(\eN)\otimes\eL^b)=0,\quad\forall t>q^\eN+q^\cL,\forall a\ges0,b\gg0.
$$
The sub-additivity property~\ref{thm:subadd-q} implies that $\eN\oplus\eL$ is $(q^\eN+q^\eL)$-ample. So the vanishing holds for $\omega_Y\otimes\eF^\vee\otimes\Sym^{a+b}(\eN\oplus\eL)$, as soon as $a+b\ges\cst{\eF}$; in particular for $a\ges 0$, $b\ges\cst{\eF}$. But the latter contains $\omega_Y\otimes\eF^\vee\otimes\Sym^{a}(\eN)\otimes\eL^b$ as a direct summand.

\noindent(ii) Since $X$ is non-singular, it is enough to consider $\eG=\omega_X\otimes\eF^\vee$, with $\eF$ locally free on $X$ (cf. \cite[Lemma III.3.2]{hart-as}). Let $\cF,\cL$ be the sheaves induced, respectively, by $\eF,\eL$ on $\hat X_Y$. The exact sequence in local cohomology for $Y\subset X$ and the formal duality (cf. \cite[Theorem III.3.3]{hart-as}) yield:
$$
\begin{array}{cl}
(1)&\,\Leftrightarrow\;
H^{\dim X-t-1}(\hat X_Y,\cF)\text{ is finite dimensional,}
\\ 
(2)&\,\Leftrightarrow\;
H^{\dim X-t-1}(\hat X_Y,\cF\otimes\cL^{-b})=0.
\end{array}
$$
It remains to apply \ref{thm:finite-dim} and (i) above, respectively.  
\end{m-proof}


\subsection{Cohomology of the complement}\label{ssct:groth}
 
In  \cite[Expos\'e XIII, Conjecture 1.3]{groth}, Grothendieck discusses the finite dimensionality of the cohomology groups of coherent sheaves on the complement of lci subvarieties in projective spaces. Hartshorne addressed the issue for smooth subvarieties (cf. \cite[Corollary 5.7]{hart-cdav}). Here we extend his result to the relative setting. 

Let $S$ be an irreducible projective variety and $\eE$ a locally free sheaf of rank $r\ges 2$ on $S$. Let $X:=\mbb P(\eE^\vee)=\Proj(\Sym^\bullet\eE)$ and $\pi:X\to S$ be the natural projection. The Euler sequence shows that the relative tangent bundle $\eT_{X,\pi}$ is $\pi$-relatively ample, so it is $\dim S$-ample. 

\begin{m-corollary}\label{cor:proj-bdl} 
Let the situation be as above. Suppose $Y\subset X$ is a smooth family of subvarieties parametrized by $S$ of relative codimension $\delta$, $\dim Y>\dim S$; that is, $\rd\pi_Y:\eT_Y\to\pi_Y^*\eT_S$ is surjective and $\codim_X(Y)=\delta$. Then $H^t(X\sm Y,\eG)$ is finite dimensional, for all coherent sheaves $\eG$ on $X\sm Y$ and $t\ges\dim S+\delta$. 
\end{m-corollary}
We remark that Hartshorne's result corresponds to $S=\{\text{point}\}$.

\begin{m-proof} 
The exact diagram 
$$\xymatrix@R=1.25em{&0\ar[d]&0\ar[d]&&\\ 0\ar[r]&\eT_{Y,\pi_Y}\ar[r]\ar[d]&\eT_Y\ar[r]\ar[d]&\pi_Y^*\eT_S\ar[r]\ar@{=}[d]&0\\ 0\ar[r]&\eT_{X,\pi}\rst_Y\ar[r]\ar[d]&\eT_X\rst_Y\ar[d]\ar[r]&\pi^*\eT_S\rst_Y\ar[r]&0\\ &\eT_{X,\pi}\rst_Y\Big/\eT_{Y,\pi_Y}\ar@{=}[r]\ar[d]&\eN_{Y/X}\ar[d]&&\\ &0&0&&}$$
shows that $\eN_{Y/X}$ is a quotient of $\eT_{X,\pi}\rst_Y$, so is $\dim S$-ample. Now apply~\ref{cor:finite-dim}. 
\end{m-proof}


\section{The G2 property}\label{sct:G2}

Here we generalize several results in \cite[Section 6]{hart-cdav}. 

\begin{m-lemma}\label{lm:L-ample}
{\rm(cf. \cite[Lemma 6.1]{hart-cdav})} 
Let $(Y,\eO_Y(1))$ be a projective variety. We consider an invertible sheaf $\eL$ and locally free sheaves $\eE,\eF$ on $Y$. Denote 
$$
\;h_{\eF}(a,b):=h^{0}(Y,\eF\otimes \Sym^{a}(\eE^{\vee})\otimes\eL^{-b}).
$$ 
The following properties are satisfied: 
\begin{enumerate}[leftmargin=5ex]
\item[\rm(i)] 
There is a polynomial $P^{\eO_Y(1),\eL,\eF}_{\dim Y+\rk\eE-1}(a,|b|)$ of total degree $\dim Y+\rk\eE-1$ such that 
$h_{\eF}(a,b)\leqslant P^{\eO_Y(1),\eL,\eF}_{\dim Y+\rk\eE-1}(a,|b|)$, for all 
$(a,b)\in\mbb Z^2$; 

\item[\rm(ii)] 
If $\eL$ is $(\dim Y-1)$-ample, then it holds: 
\begin{equation}\label{eq:b>a}
h_{\eF}(a,b)=0,\;\text{for}\;\;
b\ges\cst{\eO_Y(1),\eL,\eE}_1\cdot\,a+\cst{\eO_Y(1),\eL,\eF}_2.
\end{equation}

\item[\rm(iii)] 
If $\eE$ is $(\dim Y-1)$-ample, then it holds: 
\begin{equation}\label{eq:a>b}
h_{\eF}(a,b)=0,\;\text{for}\;\;
a\ges\cst{\eO_Y(1),\eE,\eL}_1\cdot\,b+\cst{\eO_Y(1),\eE,\eF}_2.
\end{equation}
\end{enumerate}
\end{m-lemma}

\begin{m-proof} 
We fix $\eO_Y(1)$ sufficiently ample (cf. \ref{thm:unif-q}) and consider the regularity of the various sheaves with respect to it. Also, we assume that $Y$ is irreducible, otherwise we prove the estimates on its components.
\medskip 

\nit(i) Note that $\eF^\vee\otimes\eO_Y(c_0)$ is globally generated, so $\eF$ is a subsheaf of $\eO_Y(c_0)^{\oplus\rk\eF}$, for $c_0\gg0$; actually, it suffices to take $c_0:=\max\{1,\reg\eF^\vee\}$. Thus it is enough to prove the statement for $\eF=\eO_Y(c_0)$. 
Denote $e:=\rk(\eE)$. Similarly, we have $\eE^\vee\subset\eO_Y(a_0)^{\oplus e}$ and $\eL^{\pm 1}\subset\eO_Y(b_0)$, where $a_0,b_0$ depend respectively on $\reg\eE,\reg\eL^{\pm1}$ as above.  It follows that: 
$$
\begin{array}{c}
\Sym^a(\eE^\vee)\otimes\eL^{-b}\subset
\eO_Y(aa_0)^{\oplus \bnom{a+e-1}{e-1}}\otimes\eO_Y(|b|b_0)
\\[2ex] \Rightarrow\;\;
h^0(Y,\eO_Y(c_0)\otimes\Sym^a(\eE^\vee)\otimes\eL^{-b})
\les
\bnom{a+e-1}{e-1}\cdot h^0(Y,\eO_Y(aa_0+|b|b_0+c_0)).
\end{array}
$$
It is well-known that $h^0(\eO_Y(k))$ is bounded above by a polynomial in $k$, of degree $\dim Y$ (cf. \cite[1.2.33]{laz}). But $\bnom{a+e-1}{e-1}$ is a polynomial in $a$ of degree $e-1$; overall, we obtain a polynomial in $a,|b|$ of total degree $\dim Y+e-1$. 
\medskip

\nit(ii) 
Let $\omega_{Y}$ be the dualizing sheaf of $Y$. There is $c_{0}=c_0(Y)\geqslant1$ such that $\omega_{Y}\otimes\eO_Y(c_0)$ is globally generated. A generic section induces an inclusion $\eO_Y(-c_0)\subset\omega_Y$, so:
$$
\begin{array}{rl}
h^{0}(Y,\eF\otimes\Sym^{a}(\eE^{\vee})\otimes\eL^{-b})
&
\les h^{0}(Y,\omega_{Y}\otimes\eF(c_0)\otimes\Sym^{a}(\eE^{\vee})\otimes\eL^{-b})
\\[1.5ex]&
= 
h^{\dim Y}(Y,\eF^\vee(-c_0)\otimes\Sym^{a}(\eE)\otimes\eL^b).
\end{array}
$$
\nit\textit{Claim}\quad The right hand-side vanishes for $b$ as in \eqref{eq:b>a}. 

Henceforth we replace $\eF$ by $\eF(-c_0)$ and verify the statement for $h^{\dim Y}(\eF^\vee\otimes \Sym^{a}(\eE)\otimes\eL^{b})$. In order to track the dependence of the constants on the various parameters, note that the effect of replacing  $\eF\rightsquigarrow\eF(-c_0)$ is $\reg\eF^\vee\rightsquigarrow\reg\eF^\vee-c_0$, where $c_0$ depends only on $Y$. 

We observe that it is enough to prove the claim for $Y$ reduced and for an arbitrary coherent sheaf $\eG$ on $Y$ instead of $\eF^\vee$. Indeed, for  $\eI:=\Ker(\eO_Y\to\eO_{Y_{\text{red}}})$, there is $r>0$ such that $\eI^r=0$. Thus $\eO_Y$ admits a filtration (similar to \eqref{eq:Ya}) by the quotients $\eI^{k-1}/\eI^k$, $1\les k\les r$, which are $\eO_{Y_{\text{red}}}$-modules, and we may use the estimates for $\eF^\vee\otimes(\eI^{k-1}/\eI^k)$ on $Y_{\text{red}}$, which is coherent. 

Since $Y$ is irreducible and reduced, we have $H^0(\eO_Y)=\kk$. The uniform $q$-ampleness property (cf. \ref{thm:unif-q}) implies that there are constants depending only on $\eO_Y(1)$, such that: 
$$
H^{\dim Y}(Y,\eG\otimes\Sym^a\eE\otimes\eL^b)=0,\;
\forall b\ges\cst{\eO_Y(1),\eL}_1\cdot\reg(\eG\otimes\Sym^a\eE)+\cst{\eO_Y(1),\eL}_2.
$$
Since $\Sym^a\eE$ is a direct summand of $\eE^{\otimes a}$, the sub-additivity of the regularity (cf. \ref{thm:unif-q}(ii)) yields 
$\;
\reg(\eG\otimes\Sym^a\eE)\les a\cdot\reg(\eE)+\reg(\eG),
$
so \eqref{eq:b>a} holds for: 
$$\;
b\ges\cst{\eO_Y(1),\eL}_1\cdot\,(a\cdot\reg(\eE)+\reg(\eG))+\cst{\eO_Y(1),\eL}_2.
$$ 
The coefficient of $a$ is indeed independent of $\eG$. 
\medskip 

\nit(iii) 
Again, we prove the estimate for $Y$ reduced and $\eF^\vee$ replaced by an arbitrary coherent sheaf $\eG$ on $Y$. By the uniform $(\dim Y-1)$-ampleness, $h^{\dim Y}(\eG\otimes \Sym^{a}(\eE)\otimes\eL^{b})$ vanishes for 
$$
a\ges\cst{\eO_Y(1),\eE}_1\cdot\reg(\eG\otimes\eL^b)+\cst{\eO_Y(1),\eE}_2.
$$
It remains to apply \ref{thm:unif-q}(ii): $\;\reg(\eG\otimes\eL^b)\les b\reg\eL+\reg\eG$. 
\end{m-proof}

\begin{m-proposition}\label{prop:L-any}
{\rm(cf. \cite[Theorem 6.2, Corollary 6.6]{hart-cdav})} 
Let the situation be as in \ref{XYN}. We assume that $Y$ is lci and its normal bundle is $(\dim Y-1)$-ample, of rank $\nu$. For any locally free sheaf $\cal F$ and invertible sheaf $\cal{L}$ on $\mfrak{X}$, there is a polynomial of degree $\dim Y+\nu$ such that: 
$$
h^{0}(\mfrak{X},\cal F\otimes\cal{L}^{b})\les P^{Y,\cal L,\cal F}_{\dim Y+\nu}(b),
\;\text{ for }b\gg0.
$$
\end{m-proposition}

\begin{m-proof}
Since $Y$ is lci, $\omega_Y$ is invertible. We fix a sufficiently (Koszul) ample invertible sheaf $\eA$ on $Y$, such that $\eA^{-1}\subset\omega_Y$. Let $\eF,\eL$ be the restrictions of $\cal F,\cal{L}$ to $Y$, respectively. 
With the notation of \eqref{eq:a>b}, for $\gamma:=\cst{\eA,\eN,\eL}_1+1$ and $b>\cst{\eA,\eN,\eF}_2$, we have the estimate:
$$
h^{0}(\mfrak{X},\cal F\otimes\cal{L}^{b})
\leqslant
\sum_{a=0}^{\infty}h^{0}(Y,\eF\otimes\Sym^{a}(\eN^{\vee})\otimes\eL^{b})
=
\sum_{a=0}^{\gamma b}h^{0}(Y,\eF\otimes\Sym^{a}(\eN^{\vee})\otimes\eL^{b}).
$$
Since $\eF$ is a subsheaf of $(\eA^{c_0})^{\oplus\rk\eF}$, for $c_0=\max\{1,\reg\eF^\vee\}$, it is enough to prove the proposition for $\eF=\eA^{c_0}$. 

Consider $S:=\mathbb{P}(\eO_{\mathbb{P}(\eN^{\vee})}(-1)\oplus\eO_{\mathbb{P}(\eN^{\vee})})$ and denote by $\eO_{S}(1)$ the relatively ample invertible sheaf on it. The right hand-side above can be re-written as follows:
$$
\begin{array}{l}
\disp
\text{rhs}=
\sum_{a=0}^{\gamma b} h^{0}(Y,\eA^{c_0}\otimes\Sym^{a}(\eN^{\vee})\otimes\eL^{b})
\les
\sum_{a=0}^{\gamma b} 
h^0(Y,\omega_{Y}\otimes\eA^{c_0+1}\otimes\Sym^{a}(\eN^\vee)\otimes\eL^{b})
\\[3ex]
\disp=
\sum_{a=0}^{\gamma b} 
h^{\dim Y}(Y,\eA^{-c_0-1}\otimes\Sym^{a}(\eN)\otimes\eL^{-b})
\\[3ex] \disp
=
\sum_{a=0}^{\gamma b}
h^{\dim Y}(\mathbb{P}(\eN^{\vee}),\eA^{-c_0-1}\otimes\eO_{\mathbb{P}(\eN^{\vee})}(a) \otimes\eL^{-b})
=
h^{\dim Y}(S,\eA^{-c_0-1}\otimes\eO_{S}(\gamma b)\otimes\eL^{-b}).
\end{array}
$$
But $h^{\dim Y}(S,\eO_{S}(\gamma b)\otimes\eL^{-b})$ is bounded above by a polynomial in $b$, depending on $\eO_S(\gamma)\otimes\eL^{-1}$, of degree at most $\dim S=\dim Y+\nu$ (cf. \cite[1.2.33]{laz}). In order to include the term $\eA^{-c_0-1}$, we use the exact sequence 
$$
0\to\eA^{-c_0-1}\to\eO_Y\to\eO_{Y_1}\to0,\;\;\dim Y_1=\dim Y-1,
$$
which yields: 
$$
\text{rhs}\les h^{\dim X}(\eO_S(\gamma b)\otimes\eL^{-b}) + h^{\dim X-1}(\eO_S(\gamma b)\otimes\eL^{-b}\rst_{Y_1}).
$$ 
Of course, $Y_1$ depends on $c_0$, \textit{a posteriori} on $\eF$. However, both terms are bounded above by polynomials in $b$, of degree at most $\dim Y+\nu$ and $\dim Y+\nu-1$, respectively.
\end{m-proof}

\begin{m-theorem}\label{thm:G2}
{\rm(cf. \cite[Theorem 6.7]{hart-cdav})} 
Let the situation be as in \ref{XYN}. We assume: 
\begin{itemize}[leftmargin=5ex]
\item $Y$ is connected, lci, $\dim Y\geqslant1$;
\item the normal bundle $\eN$ of $Y$ is $(\dim Y-1)$-ample.
\end{itemize}
Then the following statements hold: 
\begin{enumerate}[leftmargin=5ex]
\item[\rm(i)] 
$\trdeg_{\kk}K(\mfrak{X})\leqslant\dim Y+\rk\eN$;
\item[\rm(ii)] 
If $\trdeg_{\kk}K(\mfrak{X})=\dim Y+\rk\eN$, then $K(\mfrak{X})$ is a finitely generated extension of $\kk$.
\end{enumerate}
\end{m-theorem}

\begin{m-proof}
With our preparations, the proof is identical to \textit{loc.~cit.}. Let $\tau:=\trdeg_\kk K(\mfrak X)$ and choose a transcendence basis $\xi_1,\dots,\xi_\tau\in K(\mfrak X)$. Since $\mfrak X$ is regular, there are invertible sheaves $\cal L_i$ on $\mfrak X$, $i=1,\dots,\tau$, and sections $s_{0i},s_{1i}$ in $\cal L_i$ such that $\xi_i=s_{1i}/s_{0i}$. Then $\cal L:=\cal L_1\otimes\dots\otimes\cal L_\tau$ admits non-zero sections $s_0,s_1,\dots,s_\tau$ such that $\xi_i=s_i/s_0$. 

The ring $R:=\uset{b\ges 0}{\sum} H^0(\mfrak X,\cal L^b)$ is an integral domain, because $Y$ is connected, and one has the injective ring homomorphism: 
$$
R\to K(\mfrak X)[\xi],\quad 
\si_b\mt \frac{\si_b}{s_0^b}\cdot\xi^b,\text{ for }\si_b\in H^0(\mfrak X,\cal L^b).
$$
By Proposition~\ref{prop:L-any}, there is a polynomial $P=P^{Y,\cL}_{\dim Y+\nu}$ such that: 
$$
h^0(\mfrak X,\cal L^b)\les P(b),\text{ for }b\gg0.
$$ 
Then \cite[Lemma 6.3]{hart-cdav} implies: 
$\,\trdeg K(R)\les\dim Y+\nu+1.$ 
But $\kk[s_0,s_1,\dots,s_\tau]\subset R$, so $\tau+1\les\trdeg K(R)$, hence $\tau\les\dim Y+\nu$. 

Now assume that $\tau=\dim Y+\nu$. Then, by \textit{idem}, Lemma~6.3, 6.4, we have: 
\begin{itemize}[leftmargin=5ex]
\item[--] 
The extension $\kk(s_0,s_1,\dots,s_\tau)\hra K(R)$ is  finite.
\item[--] 
$R$ is integrally closed in $K(\mfrak X)[\xi]$ and the extension $K(R)\hra K(\mfrak X)(\xi)$ is algebraic, so $K(R)=K(\mfrak X)(\xi)$.
\end{itemize}
We conclude that $\kk(\xi_1,\dots,\xi_\tau)\hra K(\mfrak X)$ is a finite extension.
\end{m-proof}

\begin{m-corollary}\label{cor:G2}
{\rm(cf. \cite[Corollary 6.8]{hart-cdav})} 
Let $X$ be a projective scheme which is non-singular in a neighbourhood of a closed, connected, lci subscheme $Y\subset X$. We assume that the normal bundle $\eN_{Y/X}$ is $(\dim Y-1)$-ample. Then $K(\hat{X}_{Y})$ is a finite extension of $K(X)$, in other words $Y$ is G2 in $X$. 
\end{m-corollary}

\begin{m-proof}
Indeed, in this case $K(X)$ is a subfield of $K(\hat X_Y)$, so $\trdeg_\kk K(\hat X_Y)\ges\dim X=\dim Y+\nu$, hence we are in the case (ii) of the previous theorem. 
\end{m-proof}

The result is optimal, in the sense that one can not conclude that $Y$ is G3 in $X$ (cf. \cite[Example pp. 199]{hart-as}, \cite[Example 9.1]{bad}). In Section \ref{sct:G3}, we shall see several conditions which ensure the G3 property.

\subsection{A formality criterion}\label{ssct:chen}

Our discussion yields a short proof of a result obtained in \cite{chen}. One says that the \emph{formal principle} holds for a pair $(X,Y)$ consisting of a scheme $X$ and a closed subscheme $Y$ if the following condition is satisfied: for any other pair $(Z,Y)$ such that $\hat Z_Y\cong\hat X_Y$, extending the identity of $Y$, there is an isomorphism between \'etale neighbourhoods of $Y$ in $X$ and in $Z$ which induces the identity on $Y$.

\begin{m-theorem}\label{thm:chen} 
{\rm(cf. \cite[Theorem 3]{chen})} 
Let $X$ be a non-singular projective variety and $Y\subset X$ be a closed, connected, lci subscheme. Assume that $\eN_{Y/X}$ is $(\dim Y-1)$-ample. Then the formal principle holds for the pair $(X,Y)$. 
\end{m-theorem}

Note that this strengthens \textit{loc.~cit.}, since $Y$ is assumed only lci, rather than smooth. 

\begin{m-proof}
Corollary \ref{cor:G2} implies that $Y$ is G2 in $X$. But Gieseker proved (cf. \cite[Theorem 4.2]{gskr}, \cite[Corollary 9.20, 10.6]{bad}) that in this case the formal principle applies to $(X,Y)$. 
\end{m-proof}

Griffiths and Commichau-Grauert obtained similar results in complex analytic setting. On one hand, Griffiths investigates in \cite{griff} the formality/rigidity property described above, for smooth subvarieties $Y\subset X$ whose normal bundle is partially positive/negative. In [\textit{idem}, VI.~\S5, pp. 425], he considers subvarieties whose normal bundle $\eN_{Y/X}$ admits a Hermitian metric whose curvature has signature $(s,t)$, with $s+t=\dim Y$, and proves in [\textit{ibid.}, II.~\S2,\,3] the rigidity of the embedding $Y\subset X$, as soon as $s\ges 2$. 

The concept of partial positivity/negativity of the normal bundle is inspired from Andreotti-Grauert \cite{andr+grau}, who studied arbitrary vector bundles which admit Hermitian metrics whose curvature has mixed signature. Their essential cohomological property is the following: if $\eN$ admits a Hermitian metric of signature $(s,t)$, then $\eN$ is $(\dim Y-s)$-ample (cf. \cite[Proposition 28, pp.~257]{andr+grau}, \cite[(7.28), pp.~432]{griff}). 

On the other hand, Commichau-Grauert proved in \cite[Satz 4]{comm+graut} a formality/rigidity result for subvarieties with $1$-positive normal bundle. It turns out that a $1$-positive vector bundle on a smooth projective variety $Y$ is $(\dim Y-1)$-ample (cf. \cite[Satz 2]{comm+graut}). 

We conclude that the cohomological approach adopted in this article yields under weaker assumptions the rigidity results obtained in \cite{griff,comm+graut}.


\section{Examples of subvarieties with partially ample normal bundle}\label{sct:apply}

\subsection{Elementary operations}\label{ssct:elem}

\begin{m-corollary}\label{cor:elem} 
Let $X$ be a non-singular projective variety. The following statements hold:
\begin{enumerate}[leftmargin=5ex]
\item[\rm(i)] 
Let $Y_2\subset Y_1\subset X$ be connected lci, $\dim Y_2\ges1$. Suppose $\eN_{Y_2/Y_1}, \eN_{Y_1/X}$ are respectively $q_2$-, $q_1$-ample, $q_1+q_2<\dim Y_2$. Then $Y_2$ is G2 in $X$. \\ 
In particular, $Y_2\subset X$ is G2, if $\eN_{Y_1/X}$ is ample and $\eN_{Y_2/Y_1}^\vee$ is not pseudo-effective.
\smallskip

\item[\rm(ii)] 
Suppose $Y_1, Y_2$ are lci in $X$, 
$
\codim (Y_1\cap Y_2)=\codim (Y_1)+\codim(Y_2),
$
and $\eN_{Y_j/X}$ is $q_j$-ample, for $j=1,2$. Then $\eN_{Y_1\cap Y_2/X}$ is $(q_1+q_2)$-ample. 
\smallskip

\item[\rm(iii)] 
Suppose $Y_j\subset X_j$ are connected lci and $\eN_{Y_j/X_j}$ is $(\dim Y_j-1)$-ample (thus $Y_j\subset X_j$ is G2), for $j=1,2$. Then $Y_1\times Y_2$ is lci and G2 in $X_1\times X_2$. 
\smallskip

\item[\rm(iv)] 
Let $f:X'\to X$ be a surjective, flat morphism. Suppose $Y\subset X$ is lci and $\eN_{Y/X}$ is $(\dim Y-1)$-ample. Then $Y':=f^{-1}(Y)\subset X$ is lci and $\eN_{Y'/X'}$ is $(\dim Y'-1)$-ample.
\end{enumerate}
\end{m-corollary}

\begin{m-proof}
The statements are consequences of the sub-additivity properties \ref{thm:subadd-q} and \ref{lm:q-prod}. 

\nit(i) The sequence 
$0\to\eN_{Y_2/Y_1}\to\eN_{Y_2/X}\to\eN_{Y_1/X\rst_{Y_2}}\to 0$ 
implies that $\eN_{Y_2/X}$ is $(q_1+q_2)$-ample.

\nit(ii) Note that $Y_1\cap Y_2$ is lci in $X$ and its normal bundle fits into: 
$$
0\to\eN_{Y_2/X\rst_{Y_1\cap Y_2}}
\to \eN_{Y_1\cap Y_2/X}
\to \eN_{Y_1/X\rst_{Y_1\cap Y_2}}\to 0.
$$

\nit(iii) We apply Lemma~\ref{lm:q-prod} to  
$\eN_{Y_1\times Y_2/X_1\times X_2}=\eN_{Y_1/X_1}\boxplus\eN_{Y_2/X_2}$ and deduce that $\eN_{Y_1\times Y_2/X_1\times X_2}$ is $(\dim Y_1+\dim Y_2-1)$-ample. 

\nit(iv) Note that $\eN_{Y'/X'}=f^*\eN_{Y/X}$. Since $f$ is flat, its fibre dimension is constant. The claim follows from Leray's spectral sequence applied to $Y'\to Y$.
\end{m-proof}

\begin{m-corollary}\label{cor:G2-inters}
Suppose $Y_1, Y_2$ are lci in $X$, 
$\codim (Y_1\cap Y_2)=\codim (Y_1)+\codim(Y_2),$ and $\eN_{Y_j/X}$ is $q_j$-ample, for $j=1,2$. If $Y_1\cap Y_2$ is connected and $q_2<\dim(Y_1\cap Y_2)$, e.g. $q_2=0$, then $Y_1\cap Y_2$ is G2 in $Y_1$. 
\end{m-corollary}
The connectedness and G3-property of the intersection of partially ample subvarieties is discussed in Section~\ref{sct:connect}. The corollary says that the analogous G2-property holds under fairly general circumstances. 

\begin{m-proof}
Note that $\eN_{Y_1\cap Y_2/Y_1}\cong\eN_{Y_2/X\rst_{Y_1\cap Y_2}}$. 
\end{m-proof}


\subsection{Varieties whose cotangent bundle is not pseudo-effective}\label{ssct:not-pseff} 

Suppose $Y\subset X$ is a smooth subvariety. Proposition~\ref{prop:y-1} relates the $(\dim Y-1)$-ampleness of the normal bundle of $Y$ to the non-pseudo-effectiveness of $\eN_{Y/X}^\vee$. Since $\eN_{Y/X}$ is a quotient of $\eT_X$, one hopes that Theorem~\ref{thm:G2} applies to (generic) subvarieties of projective varieties $X$ such that $\eO_{\mbb P(\eT_X)}(1)$ is not pseudo-effective. Examples of projective varieties with $\eO_{\mbb P(\eT_X)}(1)$ non-pseudo-effective include uniruled varieties and, possibly, Calabi-Yau varieties. 

Indeed, the canonical bundle of an uniruled variety $X$ is not pseudo-effective (cf. \cite{bdpp}), hence its cotangent bundle is the same (cf. \cite[Theorem 6.7(a)]{dps}). One can explicitly see this as follows: it is enough to find a moving curve $\vphi=(\vphi_P,\vphi_Z):C\to\mbb P^1\times Z$, as in Proposition~\ref{prop:y-1}, on a ruling $\mbb P^1\times Z$ of $X$. Consider a complete intersection $C\subset Z$ of high degree, so $\vphi_Z$ is the inclusion, and a finite morphism $\vphi_P:C\to\mbb P^1$. Then $\vphi$ is movable because $\mbb P^1$ is homogeneous. The line bundle $\eL=\eO_C(x)$, $x\in C$, on $C$ is ample and $\eL\subset\vphi_P^*(\eT_{\mbb P^1})$, since $\vphi_P^*\eT_{\mbb P^1}$ is globally generated. Similarly, the normal bundle $\eN_{C/Z}$ is a direct sum of globally generated line bundles of large degree, so $\cHom(\eL,\vphi_Z^*\eN_{C/Z})$ is globally generated too. Hence the inclusion $\eL\subset\vphi^*\eT_{\mbb P^1\times Z}$ is movable. 

Candidates of varieties with non-pseudo-effective cotangent bundle include Calabi-Yau varieties $X$. In \cite[Corollary 6.12]{dps} it is proved that $\eO_{\mbb P(\eT_X)}(1)$ is not pseudo-effective if its non-nef locus of does not cover the whole $X$.

\begin{m-corollary}\label{cor:diag-G2}
Let $X$ be a projective variety whose cotangent bundle is not pseudo-effective. Then the diagonal $\Delta_X:=\{(x,x)\mid x\in X\}$ is G2 in the product $X\times X$. 
\end{m-corollary}
(See~\ref{cor:DX} for the G3-property of the diagonal in the quasi-homogeneous case.)
\begin{m-proof}
The normal bundle of $\Delta_X$ is isomorphic to $\eT_X$ and we conclude by~\ref{cor:G2}.
\end{m-proof}

\begin{m-notation}
For shorthand, denote $\mbb P:=\mbb P(\eT_X)$ and $\mbb P_Y:=\mbb P(\eT_{X\rst_Y})$ its restriction to $Y$; let $\pi:\mbb P\to X$ be the projection. 
Define ${\rm Mov}(\mbb P_Y)_{\mbb Q}\subset H_2(\mbb P_Y;\mbb Q)$ to be the cone generated by the classes of movable curves on $\mbb P_Y$ and ${\rm Mov}(\mbb P)_{\mbb Q}$ similarly.
\end{m-notation}

If $\eO_{\mbb P}(1)$ is not pseudo-effective, there are numerous subvarieties $Y\subset X$ with $\eO_{\mbb P_Y}(1)$ non-pseudo-effective. Indeed, there is a reduced, movable curve $M\subset\mbb P$ such that $\eO_{\mbb P}(1)\cdot M<0$, of the form $M=\si_*\tld M$, with $\tld{\mbb P}\srel{\si}{\to}\mbb P$ birational and $\tld M$ a complete intersection on $\tld{\mbb P}$. Take $Y\subset X$ such that $M\subset\mbb P_Y$. Then $\tld M\subset\tld{\mbb P}_Y:=\mbb P_Y\times_{\mbb P}\tld{\mbb P}\subset\tld{\mbb P}$ is a complete intersection, so $M=\si_*\tld M$ is movable on $\mbb P_Y$ and $\eO_{\mbb P_Y}(1)\cdot M<0$. 

\begin{m-theorem}\label{thm:mov}
Let the notation be as above and assume that $\eO_{\mbb P}(1)$ is not pseudo-effective. Consider a smooth irreducible subvariety $Y\subset X$ such that $\eO_{\mbb P_Y}(1)$ is not pseudo-effective, too. Then $\eN_{Y/X}$ is $(dim Y-1)$-ample, so $Y$ is G2.
\end{m-theorem}

\begin{m-proof}
Indeed, $\eT_{X\rst_Y}$ is $(\dim Y-1)$-ample, so $\eN_{Y/X}$ is the same.
\end{m-proof}

One may ask if one can weaken the condition on $\eO_{\mbb P}(1)$. The answer is negative, for sufficiently general subvarieties. Hence, if $X$ has pseudo-effective cotangent bundle and $Y\subset X$ is a general subvariety, the partial ampleness of $\eN_{Y/X}$ has to be verified directly. 

\begin{m-proposition}\label{prop:mov} 
Let $Y\subset X$ be a positive-dimensional subvariety and 
$$\iota_*:H_2(\mbb P_Y;\mbb Q)\to H_2(\mbb P;\mbb Q)$$ 
be the homomorphism induced by the inclusion. In the situations enumerated below, holds: 
\begin{equation}\label{eq:mov}
\iota_*\big({\rm Mov}(\mbb P_Y)\big)_{\mbb Q}\subseteq{\rm Mov}(\mbb P)_{\mbb Q}.
\end{equation}
\begin{enumerate}[leftmargin=5ex]
\item[\rm(i)] 
An algebraic group $G$ acts on $X$ with an open orbit $O$, such that the stabilizer of a point $x\in O$ acts with open orbit on $\eT_{X,x}$ and $Y\cap O\neq\emptyset$. 
\item[\rm(ii)] 
$Y$ is a movable, very general subvariety of $X$ and $\kk$ is uncountable. 
\end{enumerate}
Hence, if $\eO_{\mbb P_Y}(1)$ is not pseudo-effective, then $\eO_{\mbb P}(1)$ is the same.
\end{m-proposition}

\begin{m-proof} 
(i) The hypothesis implies that the $G$-translates of a movable curve on $\mbb P_Y$ cover an open subset of $\mbb P(\eT_X)$. 

\nit(ii) Let $S$ be an affine parameter space, $Y_S\subset S\times X$ an $S$-flat family of subvarieties of $X$. Curves on $X$ are parametrized by their Hilbert polynomial with respect to $\eO_X(1)$, of degree one, with integer coefficients. Let $\Pi$ be the set of polynomials which occur as Hilbert polynomials of movable curves on $Y_s$, for $s\in S$; it is a countable set. 

For $P\in\Pi$, denote by $\Hilb^P_{Y_S/S}\srel{\pi_P\,}{\to}S$ the corresponding relative Hilbert scheme; it is projective over $S$. We are interested only in the components corresponding to curves. For $s\in S$, let $\Pi_s\subset \Pi$ be the set of polynomials $P_s$ such that $\pi_{P_s}$ is not dominant; denote $\Pi_{\text{rigid}}:=\uset{s\in S}{\bigcup}\Pi_s$. The image of $\Hilb^{\Pi_{\text{rigid}}}_{Y_S/S}\to S$ is a countable union of proper subvarieties. 

Take $s'\in S$ in the complement ($\kk$ is uncountable); let $P_{s'}$ be the Hilbert polynomial of some movable curve $C_{s'}\subset Y_{s'}$. Then $P_{s'}\not\in\Pi_{\text{rigid}}$, by the choice of $s'$, so $\Hilb^{P_{s'}}_{Y_S/S}\srel{\pi_{P_{s'}}\;}{\lar}S$ dominates $S$, hence $\pi_{P_{s'}}$ is surjective. Let $\Pi':=\Pi\sm\Pi_{\text{rigid}}$. The components of $\Hilb^{\Pi'}_{Y_S/S}$ (corresponding to movable curves) dominate $S$, so they are flat over the very general point $o\in S$. 

We claim that movable curves on $Y_o$ are movable on $X$. Indeed, for $P_o$ as above, consider the universal curve $\cal C_S\subset\Hilb^{P_o}_{Y_S/S}\times_S Y_S$. The family $\cal C_{o}\subset \Hilb^{P_o}_{Y_o}\times Y_o$ covers an open subset of $Y_o$. The continuity of the $S$-morphism $\Hilb^{P_o}_{Y_S/S}\times_S Y_S\to Y_S$ implies that the same holds for $\cal C_{s}\subset \Hilb^{P_o}_{Y_s}\times Y_s$, for $s$ in a neighbourhood of $o\in S$. Finally, $Y_S\to X$ is dominant, so $\cal C_S$ covers an open subset of $X$. 
\end{m-proof}

For $\kk=\mbb C$ and $X$ K\"ahlerian, results obtained by Boucksom imply: if $\eO_{\mbb P}(1)$ is pseudo-effective and $Y\subset X$ is such that $\mbb P_Y$ is not contained in the non-nef locus of $\eO_{\mbb P}(1)$, then $\eO_{\mbb P_Y}(1)$ is pseudo-effective too. For the converse, however, it is not clear on which subvarieties should one check the non-pseudo-effectiveness of $\eO_{\mbb P}(1)$.


\section*{\textbf{Part~II:\hspace{2ex}\textit{q}-ample subvarieties of projective varieties}}

Ottem introduced in \cite{ottm} the notion of an ample subvariety of a projective variety and studied the corresponding properties. In what follows, we generalize this concept and define partial ampleness for subvarieties. A number of results extend to this setting. 

Related to the discussion in Part~I, the normal bundle of partially ample, lci subvarieties are partially ample, but, additionally, one controls the cohomological dimension of their complement. Therefore, partially ample subvarieties enjoy the stronger G3 property. 

\section{Definition and first properties}\label{sct:1-property}

\begin{m-definition}\label{def:q1}{\rm(cf. \cite[Definition 3.1]{ottm})} 
Let $X$ be a projective variety over the ground field $\kk$ and $Y\subset X$ a subscheme of codimension $\delta$. We denote $\tld X:=\Bl_Y(X)$ the blow-up of $\eI_Y$ and $E_Y\subset\tld X$ the exceptional divisor. 

We say that $Y$ is a \emph{$q$-ample subscheme} of $X$ if the invertible sheaf $\eO_{\tld X}(E_Y)$ is $(q+\delta-1)$-ample. 
That is, for any locally free, hence for any coherent, sheaf $\tld\eF$ on $\tld X$ holds: 
\begin{equation}\label{eq:q11}
H^{t}(\tld X,\tld\eF\otimes\eO_{\tld X}(mE_Y))=0,\;\forall\,t\ges q+\delta,\;\forall\,m\gg0.
\end{equation}
\end{m-definition} 

\begin{m-remark}
\nit{\rm(i)} For $q=0$ one recovers the ample subschemes introduced by Ottem. 

\nit{\rm(ii)} Ample subvarieties are equidimensional (cf. \cite[Proposition 3.4]{ottm}). This is not necessarily true for $q>0$. 
Indeed, let $X:=\mbb P^2$ and $Y:=\{x=0\}\cup\{y=z=0\}$. Then $\tld X=\Bl_Y(\mbb P^2)$ is isomorphic to the blow-up of $\mbb P^2$ at $[1:0:0]$---we denote it by $\wtld{\mbb P^2}$ and by $E$ the exceptional divisor---, and $\eO_{\tld X}(E_Y)=\eO_{\mbb P^2}(1)\otimes\eO_{\wtld{\mbb P^2}}(E)$. For $m\ges 1$, the exact sequence 
$$
0\to H^1\big(\eO_{\mbb P^2}(m)\otimes\eO_{\wtld{\mbb P^2}}((m-1)E)\big)\to H^1\big(\eO_{\mbb P^2}(m)\otimes\eO_{\wtld{\mbb P^2}}(mE)\big)\to H^1(\eO_E(-m))\to 0,
$$
shows that the middle term does not vanish, so $\eO_{\tld X}(E_Y)$ is $1$-ample, hence $Y\subset\mbb P^2$ is $1$-ample. 
\end{m-remark}

\begin{m-proposition}\label{prop:flat+reduct} 
\nit{\rm(i)} 
Let $f : \cX\to S$ be flat, projective morphism and let $\cY$ be a lci closed subscheme of X such that $f_\cY :\cY \to S$ is flat. Assume that there is a point $o\in S$ such that $Y_{o}$ is $q$-ample in $X_{o}$. Then there is an open neighbourhood $U$ of $o$ such that for each $s\in U$, $Y_s$ is $q$-ample in $X_s$. 

\nit{\rm(ii)} 
Let $Y$ be a subscheme of a projective variety $X$ and let $\bar Y$ be its integral closure. Then $Y$ is $q$-ample in $X$ if and only if $\,\bar Y$ is. In particular, a subscheme $Y$ is $q$-ample if and only if the subscheme $Y$ associated to a reduction is.
\end{m-proposition}

\begin{m-proof}
See \cite{ottm}, Theorem~6.1 and Proposition~6.8, respectively. 
\end{m-proof}

\begin{m-proposition}{\rm(cf. \cite[Theorem 5.4]{ottm})}\label{prop:EY}
A subscheme $Y\subset X$ is $q$-ample if and only if the following conditions are satisfied:
$$
\biggl\{\begin{array}{l}
\eO_{E_Y}(E_Y)\text{ is $(\delta+q-1)$-ample},
\\[1ex]
\cd(X\sm Y)\les\delta+q-1.
\end{array}\biggr.
$$
(Recall that $\cd(X\sm Y)\ges\delta-1$ for any $Y$ of codimension $\delta$.)
\end{m-proposition}

\begin{m-proof}  
\nit$(\Rightarrow)$ 
By \cite[Theorem 6.3]{tot}, it is enough to check the partial ampleness property for sheaves $\tld\eF_{E_Y}$, where $\tld\eF={\tld\eA}^{-k}$, $k\ges 1$, $\tld\eA\in\Pic(\tld X)$ is ample. The sequence 
\begin{equation}\label{eq:E}
0\to\tld\eF((m-1)E_Y)\to\tld\eF(mE_Y)\to\tld\eF_{E_Y}(mE_Y)\to0,
\end{equation}
implies $H^t(\tld\eF_{E_Y}(mE_Y))=0$ for $t\ges\delta+q$ and $m\gg0$. Second, \cite[(5.1)]{ottm} implies that $H^t(X\sm Y,\tld\eF)=\disp\varinjlim H^t(\tld X,\tld\eF(mE_Y))$, which vanishes for $t\ges\delta+q$. 

\nit$(\Leftarrow)$ Let $\tld\eF$ be a locally free sheaf on $\tld X$. For $t\ges\delta+q$, \eqref{eq:E} shows that the dimension of $H^t(\tld X,\tld\eF(mE_Y))$ is eventually constant. But the limit is $H^t(X\sm Y,\tld\eF)$ which vanishes, by the assumption on the cohomological dimension. 
\end{m-proof}

The proposition breaks the issue of deciding the partial ampleness of a subscheme $Y\subset X$ into a local and a global problem.  The partial ampleness of the normal sheaf is typically easier to verify. Unfortunately, it is more difficult to control the cohomological dimension of the complement (cf. \cite{ogus,lyub-licht}). 

For $\kk=\bbC$, a result due to Andreotti-Grauert \cite[Corollaire, pp. 250]{andr+grau} says that $\cd(X\sm Y)<q$ if $X\sm Y$ is a strongly $q$-complete analytic variety. (That is, $X\sm Y$ admits an exhaustion function which is strongly $q$-convex.)

\begin{m-proposition}\label{prop:connect}
Assume that $X$ is a Cohen-Macaulay variety and $Y$ is $(\dim Y-1)$-ample. Then the following statements hold:
\begin{enumerate}[leftmargin=5ex]
\item[\rm(i)]
$Y$ is connected;
\item[\rm(ii)]
If moreover $Y$ is Cohen-Macaulay, then $Y$ is also equidimensional.
\end{enumerate}
\end{m-proposition}

\begin{m-proof}
(i) Proposition~\ref{prop:EY} implies that $\cd(X\sm Y)\les\dim X-2$. According to \cite[Ch.~III, Theorem~3.4]{hart-as}, $H^0(X,\eO_X)\to H^0(\hat X_Y,\cO_{\hat X_Y})$ is an isomorphism, so $Y$ is connected. 
(Note that the reference requires $X$ to be smooth. However, this assumption is used only to apply the Serre duality, cf. \textit{loc. cit.}, proof of Theorem 3.3.)

\nit(ii) By the unmixedness theorem, local Cohen-Macaulay rings are equidimensional.
\end{m-proof}

\begin{m-proposition}\label{prop:q1}
Let $Y\subset X$ be a subscheme. We consider the conditions: 
\begin{enumerate}[leftmargin=5ex]
\item[\rm(a)]
$Y$ is $q$-ample;
\item[\rm(b)]
For all locally free sheaves $\eF$ on $X$, we have: 
\begin{equation}\label{eq:q12}
H^t(X,\eF\otimes\eI_Y^m)=0,\;\forall\,t\les \dim Y- q,\;\forall\,m\ges\cst{\eF}.
\end{equation}
\end{enumerate}
Then the following statements hold:
\begin{enumerate}[leftmargin=5ex]
\item[\rm(i)]
If $\tld X$ is Cohen-Macaulay, then $\;\text{\rm(a)}\,\Rightarrow\,\text{\rm(b)}$.\\ 
In particular, if $q\les\dim Y-1$, $Y$ is connected.
\item[\rm(ii)]
If $\tld X$ is Gorenstein, then $\;\text{\rm(a)}\,\Leftrightarrow\,\text{\rm(b)}$. 
\end{enumerate}
\end{m-proposition}

\begin{proof}
(i) Since $\eO_{\tld X}(-E_Y)$ is relatively ample for $\tld X\to X$, for $m\gg0$, Leray's spectral sequence and the Serre duality on $\tld X$ yield the following isomorphisms:
$$
H^t(X,\eF\otimes\eI_Y^m)
\cong H^t(\tld X,\eF\otimes\eO_{\tld X}(-mE_Y))
\cong H^{\dim X-t}(\tld X,\omega_{\tld X}\otimes\eF^\vee\otimes\eO_{\tld X}(mE_Y)).
$$

\nit(ii) In this case $\omega_{\tld X}$ is an invertible sheaf. The previous equation shows that the condition \eqref{eq:q11} holds for invertible sheaves $\omega_{\tld X}\otimes\eL$, with $\eL\in\Pic(X)$; we need to prove that it holds for an arbitrary coherent sheaf $\tld\eF$ on $\tld X$. 

We consider $\eA\in\Pic(X)$ ample such that $\eA(-E_Y)$ is ample on $\tld X$. For $c>0$ such that $(\tld\eF\otimes\omega_{\tld X}^{-1})\otimes\eA(-E_Y)^c$ is globally generated, we have the exact sequence:
$$
0\to\tld\eF_1:=\Ker(\veps)\to
\big(\omega_{\tld X}\otimes\eA^{-c}\otimes\eO_{\tld X}(cE_Y)\big)^{\oplus N}
\srel{\veps}{\to}\tld\eF\to0,\;\; (\text{for some }N>0).
$$
The inductive argument in \cite[Lemma~2.1]{ottm} yields the conclusion. Indeed, the previous discussion shows that $H^j(\tld\eF(mE_Y))\subset H^{j+1}(\tld\eF_1(mE_Y))$, for $j\ges\codim Y+q$ and $m\gg0$. By repeating the argument for $\tld\eF_1$, \textit{etc.}, we obtain the desired cohomology vanishing for $\tld\eF$.
\end{proof}

\begin{m-remark}
(i) In conjunction with \ref{prop:q1}(i), recall that any projective scheme $X$ admits a birational modification $\tld X\to X$ such that $\tld X$ is Cohen-Macaulay (cf. \cite{kaw}). This can be interpreted as the blow-up of some subscheme $Y\subset X$.

\nit(ii) The Gorenstein property of the blow-up and of the Rees algebra has been investigated by several authors (cf. \cite{hyry} and the references therein). An important situation, which covers many geometric applications, is when $X$ is smooth and $Y\subset X$ is lci. 
\end{m-remark}

As we shall see, sometimes it is convenient to work with the `positivity', as opposed to the `ampleness', of various objects. For this reason, we introduce the following \textit{ad~hoc} terminology. 

\begin{m-definition}\label{def:p>0}
The subscheme $Y$ of $X$ \emph{is (has the property) $p^\pos$} if, for all locally free sheaves $\eF$ on $X$, holds: 
$\;H^t(X,\eF\otimes\eI_Y^m)=0,\;\forall\,t\les p,\;\forall\,m\ges\cst{\eF}.$
\end{m-definition}
If $\tld X$ is Gorenstein, our discussion shows:  
$\;\text{$Y$ is $p^\pos\;\Leftrightarrow\;Y$ is $(\dim Y-p)$-ample}.$


\subsection*{A simple criterion}\label{ssct:criter} 

The computation of the ampleness of a subvariety is not straightforward, in general. However, this task becomes easy in the situation described below. In the sections \ref{ssct:glob-gen}, \ref{ssct:fixed} we apply the criterion to zero loci of sections in globally generated vector bundles and to sources of $G_m$-actions, respectively. 

\begin{m-proposition}\label{prop:x-b}
Let $Y$ be a subscheme of $X$. Assume that there is a scheme $V$, and a morphism 
$$\phi:\tld X=\Bl_Y(X)\to V$$ 
such that $\eO_{\tld X}(E_Y)$ is $\phi$-relatively ample. Then $Y\subset X$ is $q$-ample, for 
$$
q:=1+\dim\phi(\tld X)-\codim_X(Y).
$$
If $\tld X$ is Gorenstein (e.g. $X$ is smooth, $Y$ is lci), then $Y$ is $p^\pos$, for $p:=\dim X-\dim\phi(\tld X)-1$. 
\end{m-proposition}

\begin{proof}
Let $\tld\eF$ be a coherent sheaf on $\tld X$. Since $\eO_{\tld X}(E_Y)$ is relatively ample, it holds: 
$$
R^t\phi_*(\tld\eF\otimes\eO_{\tld X}(mE_Y))=0,\quad t>0,\;\;m\gg0.
$$
For $j\ges\codim_X(Y)+q>\dim\phi(\tld X)$, we have 
$$
H^{j}\big(\tld X,\tld\eF\otimes\eO_{\tld X}(mE_Y)\big)=
H^j\big(\,V,\phi_*(\tld\eF\otimes\eO_{\tld X}(mE_Y))\,\big).
$$ 
But the right hand-side vanishes, because $\Supp \phi_*(\tld\eF\otimes\eO_{\tld X}(mE_Y))$ is at most $\dim\phi(\tld X)$-dimensional. 
\end{proof}

\begin{m-proposition}\label{prop:bott}
Let the situation be as in~\ref{prop:x-b}, $\kk=\bbC$, and $X,Y$ smooth. Then holds: 
$$
H^t(X;\mbb Z)\to H^t(Y;\mbb Z)\text{ is: }
\biggl\{\begin{array}{rl}
\text{an isomorphism, for}&t\les p-1;
\\[1.25ex] 
\text{surjective, for}&t=p.
\end{array}
$$
In particular, if $p\ges 3$, then $\Pic(X)\to\Pic(Y)$ is an isomorphism. 
\end{m-proposition}

\begin{m-proof}
We claim that $\eO_{\tld X}(E_Y)$ is $\dim\phi(\tld X)$-positive. Indeed, $\eO_{\tld X}(E_Y)$ is $\phi$-relatively ample, so there is an embedding $\tld X\srel{\iota}{\to}\mbb P^N\times V$ (over $V$) such that $\eO_{\tld X}(m_0E_Y)=\iota^*(\eO_{\mbb P^N}(1)\boxtimes\euf M)$, 
for some $m_0>0$, $\euf M\in\Pic(V)$. Take the Fubini-Study metric on $\eO_{\mbb P^N}(1)$ and an arbitrary one on $\euf M$. We conclude by Theorem~\ref{thm:mats-bott}(ii). 
\end{m-proof}


\section{Equivalent characterization and elementary operations}\label{sct:N+cd}

The goal of this section is to study the behaviour of partial ampleness under various natural operations: intersection, pull-back, product. Combined with Theorem~\ref{thm:G3}, one obtains various sufficient criteria for a subvariety to be G3 in the ambient space.

\begin{m-asmp}
In this section, we assume that $X$ is smooth and $Y$ is lci. 
\end{m-asmp}

\begin{m-proposition}\label{prop:pull-back}
Let $f:X'\to X$ be a flat, surjective morphism, with $X,X'$ smooth. If $Y\subset X$ is lci and $p^\pos$, then $Y':=f^{-1}(Y)\subset X'$ is the same. 
\end{m-proposition}

\begin{m-proof}
Since $f$ is flat, $Y'$ is lci in $X'$ and $\codim_{X'}(Y')=\codim_{X}(Y)=\delta$. We check the property \eqref{eq:q11}. The universality property of the blow-up (cf. \cite[Ch. II, Corollary 7.15]{hart-ag})  yields the commutative diagram: 
$$
\xymatrix@R=1.5em@C=4em{
\tld X'=\Bl_{Y'}(X')\ar[d]\ar[r]^-{\tld f}&\tld X=\Bl_{Y}(X)\ar[d]
\\
X'\ar[r]^-f&X.
}
$$
The fibres of $f$ and $\tld f$ have the same dimension---let it be $d$---and $\tld f^*\eO_{\tld X}(E_{Y})=\eO_{\tld X'}(E_{Y'})$. For any coherent sheaf $\tld\eG$ on $\tld X'$ holds: 
$$
R^i\tld f_*(\tld\eG\otimes\eO_{\tld X'}(mE_{Y'}))
=R^i f_*\tld\eG\otimes\eO_{\tld X}(mE_{Y}),\;\; R^{j}\tld f_*\tld\eG=0,\,j>d.
$$ 
As $Y\subset X$ is $p^\pos$, we deduce: 
\\[1ex] \centerline{
$H^t\big(\tld X,R^i\tld f_*\tld\eG\otimes\eO_{\tld X}(mE_{Y})\big)=0$, for $i=0,\ldots,d$, $t\ges(\dim X-p)$, and $m\gg0$. 
}\\[1ex]
The Leray spectral sequence implies that $E_{Y'}$ is $\big((\dim X-p)+d-1\big)$-ample. 
\end{m-proof}

\begin{m-proposition}\label{prop:N+cd}
For $Y\subset X$ irreducible and lci, it holds: 
\begin{equation}\label{eq:equiv-p}
\text{$Y$ is $p^\pos$}
\quad\Leftrightarrow\quad 
\left\{\begin{array}{l}
\text{the normal sheaf $\eN=\eN_{Y/X}$ is $(\dim Y-p)$-ample,} 
\\[1.5ex] 
\text{the cohomological dimension $\cd(X\sm Y)\les\dim X-(p+1)$}. 
\end{array}\right.
\end{equation}
\end{m-proposition}

\begin{m-proof}
Apply~\ref{prop:EY}: 
$E_Y=\Proj\big(\Sym^\bullet(\underbrace{\eI_Y/\eI_Y^2}_{=\eN^\vee})\big)=\mbb P(\eN)$ 
and 
$\eO_{E_Y}(E_Y)=\eO_{\mbb P(\eN)}(-1).$ 
\end{m-proof}

\begin{m-proposition}\label{prop:trans}
Suppose $Z\subset Y$ is $p^\pos$, $Y\subset X$ is $r^\pos$, and both are irreducible lci. Then it holds: 
\begin{equation}\label{eq:trans-p}
\left\{
\begin{array}{l}
\text{$\eN_{Z/X}$ is $\big(\dim Y+\dim Z-(r+p)\big)$-ample,}
\\[1.5ex] 
\text{$\cd(X\sm Z)\les\dim X-(\min\{r,p\}+1)$}.
\end{array}\right.
\end{equation}
In particular,  $Z\subset X$ is $\big(p-(\dim Y-r)\big)^\pos$. 
\end{m-proposition}

\begin{m-proof}
The first inequality is proved in \ref{cor:elem}(i). The bound on the cohomological dimension is analogous to \cite[Proposition 6.4]{ottm}; we recall the proof here. Let $U_Z:=X\sm Z, U_Y:=X\sm Y$ and consider an arbitrary sheaf $\eG$ on $X$. In the exact sequence 
$$
\ldots\to H^i_{Y\sm Z}(U_Z,\eG)\to H^i(U_Z,\eG)\to H^i(U_Y,\eG)\to\ldots,
$$
the right hand-side vanishes for $i\ges\dim X-r$, because $Y\subset X$ is $r^\pos$. We claim that the left hand-side vanishes too, for $i\ges \dim X-p$. Indeed, it can be computed by using the spectral sequence (cf. \cite[Expos\'e I, Th\'eor\`eme 2.6]{groth}): 
$$
H^b(U_Z,\cal H^a_{Y\sm Z}(\eG))\;\Rightarrow\;H^{a+b}_{Y\sm Z}(U_Z,\eG),
$$
where $\cal H^a_{Y\sm Z}(\eG)$ stands for the local cohomology sheaf with support on $Y\sm Z$. The term on the left has the following properties: 
\begin{itemize}[leftmargin=5ex]
\item[$\bullet$] 
$\cal H^a_{Y\sm Z}(\eG)=\underset{m}{\varinjlim}\,
{\cE}{\kern-1pt}xt^a(\eO_{U_Z}/\eI_{Y\sm Z}^m,\eG)$ 
(cf. \cite[Expos\'e II, Th\'eor\`eme 2]{groth}), the ${\cE}{\kern-1pt}xt$ groups are supported on $Y\sm Z$, and $Z\subset Y$ is $p^\pos$, hence $H^b(U_Z,\cal H^a_{Y\sm Z}(\eG))=0,\;\forall\,b\ges\dim Y-p;$
\item[$\bullet$]
$\cal H^a_{Y\sm Z}(\eG)=0,\;\forall\,a\ges\dim X-\dim Y+1,$ because $Y\subset X$ is lci. 
\end{itemize}

We deduce: $H^i_{Y\sm Z}(U_Z,\eG)=0$, for $i\ges(\dim X-\dim Y)+(\dim Y-p-1)+1$. 
\end{m-proof}

The lack of sufficient positivity of the normal bundle $\eN_{Z/X}$ prevented us to conclude that $Z\subset X$ is $\min\{r,p\}^\pos$. However, we shall see in Proposition~\ref{prop:approx-p} that this is close to be true.

\begin{m-proposition}\label{prop:product}
\begin{enumerate}
\item[\rm(i)] 
Let $Y_1,Y_2\subset X$ be respectively $q_1$-, $q_2$-ample lci subvarieties such that 
$\;\codim(Y_1\cap Y_2)=\codim(Y_1)+\codim(Y_2).$ 
Then $Y_1\cap Y_2\subset X$ is $(q_1+q_2)$-ample. 
\item[\rm(ii)]
Suppose $Y_j\subset X_j$ are lci and $p_j^\pos$, for $j=1,2$. Then $Y_1\times Y_2\subset X_1\times X_2$ is $\min\{p_1,p_2\}^\pos$.
\end{enumerate}
\end{m-proposition}

\begin{m-proof}
(i) By \ref{cor:elem}, $\eN_{Y_1\cap Y_2/X}$ is $(q_1+q_2)$-ample. The Mayer-Vietoris sequence for a coherent sheaf $\eG$ on $X\sm (Y_1\cap Y_2)$ is: 
$$
\to H^{i-1}\big(X\sm (Y_1\cup Y_2),\eG\big)\to H^i\big(X\sm (Y_1\cap Y_2),\eG\big)\to H^i\big(X\sm Y_1,\eG\big)\oplus H^i(X\sm Y_2,\eG)\to\;.
$$
Since $X\sm(Y_1\cup Y_2)$ is closed in $(X\sm Y_1)\times(X\sm Y_2)$, we have:
$$
\cd\big(X\sm (Y_1\cup Y_2)\big)\les\cd(X\sm Y_1)+\cd(X\sm Y_2)
\les q_1+q_2+\delta_1+\delta_2-2.
$$ 
It follows $\cd\big(X\sm(Y_1\cap Y_2)\big)\les q_1+q_2+\codim (Y_1\cap Y_2)-1$, hence $Y_1\cap Y_2$ is $(q_1+q_2)$-ample (cf. Proposition \ref{prop:N+cd}). 

\nit(ii) We showed in \ref{cor:elem} that $\eN_{Y_1\times Y_2/X_1\times X_2}$ is $q$-ample, with  $q=\dim Y_1+\dim Y_2-\min\{p_1,p_2\}$. Second, we have: 
$\;X_1\times X_2\sm Y_1\times Y_2=\underbrace{(X_1\sm Y_1)\times X_2}_{=O_1}
\;\cup\;\underbrace{X_1\times(X_2\sm Y_2)}_{=O_2}.$

Since $\cd(O_1\cap O_2)<\dim(X_1\times X_2)-(p_1+p_2)-1$, $\cd(O_j)<\dim(X_1\times X_2)-p_j$, $j=1,2$, the Mayer-Vietoris sequence implies that $\cd(O_1\cup O_2)<\dim(X_1\times X_2)-\min\{p_1,p_2\}$. 
\end{m-proof}


\section{Weak positivity}\label{sct:aprox-q}

Throughout this section we assume that $X$ is a non-singular projective variety. We define a weak positivity property for a subscheme, suggested by the condition \ref{def:p>0}. This concept allows to prove a sort of transitivity for the $p^\pos$-property (cf. Proposition~\ref{prop:approx-p}). 

\begin{m-definition}\label{def:ap>0} 
We say that a subscheme $Y\subset X$ is $p^\apos$ (\emph{weakly} $p^\pos$) if there is a decreasing sequence of sheaves of ideals $\{\eJ_m\}_{m}$ such that the following conditions hold:  
\begin{equation}\label{eq:ap}
\begin{array}{lll}
\bullet\;&
\forall\,m,n\ges1\;\exists\,m'>m,\,n'>n\text{ such that }
\eJ_{m'}\subset\eI_Y^m,\;\eI_Y^{n'}\subset\eJ_{n};
&
\\[1ex] 
\bullet\;&
\text{for any locally free sheaf $\eF$ on $X$,}
\\[.5ex]
&
\exists\,\cst{\eF}\ges1\text{ such that }H^t\big(X,\eF\otimes{\eJ_{m}}\big)=0,
\;\forall\,t\les p\;\forall\,m\ges\cst{\eF}.
&\kern10ex\null
\end{array}
\end{equation}
Obviously, $Y$ is $p^\apos$ if and only if $Y_{\text{red}}$ is $p^\apos$. 
\end{m-definition}

\begin{m-lemma}\label{lm:cd-apos}
Assume that $Y\subset X$ is lci and $p^\apos$. Then we have: 
$$\cd(X\sm Y)\les\dim X-(p+1).$$
\end{m-lemma}

\begin{m-proof}
Since $X$ is non-singular, it holds (cf. \cite[Proposition III.3.1]{hart-as}): 
$$
\cd(X\sm Y)< c\;\Leftrightarrow\; 
H^t(X\sm Y,\eL)=0,\,\forall \eL\in\Pic(X),\,\forall t\ges c. 
$$
We have $X\sm Y\cong\tld X\sm E_Y$ and 
$H^t(\tld X\sm E_Y,\eL)=\varinjlim H^t(\tld X,\eL(mE_Y)),$ (cf. \cite[(5.1)]{ottm}). Since  $\Pic(\tld X)\cong\Pic(X)\oplus\mbb ZE_Y$, we have $\omega_{\tld X}\otimes\eL^{-1}\cong \eM(lE_Y)$ for some $\eM\in\Pic(X)$, $l\in\mbb Z$. The Serre duality on $\tld X$ implies that it is enough to verify:
$$
\varprojlim H^{j}(X,\eM\otimes\eI_Y^m)=0,\;\forall\eM\in\Pic(X),\;\forall j\les\dim X-c.
$$
The defining property of the sequence $\{\eJ_n\}_n$ together with the universality property of the projective limit yields: 
$\,\varprojlim H^j(X,\eM\otimes\eI_Y^m)=\varprojlim H^j(X,\eM\otimes\eJ_n).$
By hypothesis, the right hand-side vanishes for $j\les p$. 
\end{m-proof}

\begin{m-corollary}
Let $X$ be a non-singular projective variety over $\bbC$ and $Y$ a non-singular $p^\apos$ subvariety. 

\nit{\rm(i)} (cf. \cite[Corollary 5.2]{ottm}) 
$H^t(X;\mbb Q)\to H^t(Y;\mbb Q)$ is an isomorphism, for $t\les p-1$, and injective for $t=p$. 

\nit{\rm(ii)} Suppose $p\ges 3$. Then the following statements hold:

\begin{itemize}
\item[\rm(a)] The homomorphism $\Pic(X)\to\Pic(\hat{X}_Y)$ is bijective. 

\item[\rm(b)] $\Pic^0(X)\to\Pic^0(Y)$ is a finite morphism and $\NS(X)\to\NS(Y)$ has finite index.\\ 
Hence, $\Pic^0(X)=0$ implies $\Pic^0(Y)=0$ and $\Pic(Y)$ is finitely generated {\rm(cf. \cite[Corollary 8.6]{hart-cdav})}. 
\end{itemize}
\end{m-corollary}

\begin{m-proof}
(i) See \textit{loc.~cit.}. 

\nit(ii)(a) By \cite[Lemma 8.3]{hart-cdav}, one has the exact sequence: 
$$
H^1(X;\mbb Z)\to H^1(\hat X_Y;\cO_{\hat X_Y})\to\Pic(\hat X_Y)\to H^2(X;\mbb Z)\to H^2(\hat X_Y;\cO_{\hat X_Y}).
$$
But $\disp H^j(\hat X_Y;\cO_{\hat X_Y})=\varprojlim_m H^j(\eO_X/\eI_Y^m) 
=\varprojlim_n H^j(\eO_X/\eJ_n)=H^j(\eO_X)$, for $j=1,2$. Note that $\Pic(X)$ fits into a similar sequence. The previous step and the five lemma yields the conclusion. 

\nit(ii)(b) The Hodge decomposition implies that $H^t(X,\eO_X)\to H^t(Y,\eO_Y)$ are isomorphisms, for $t=1,2$. It remains to use the exponential sequences for $X$ and $Y$. 
\end{m-proof}

\begin{m-proposition}\label{prop:approx-p}
Let $Z\subset Y$ and $Y\subset X$ be non-singular $p^\pos$ subvarieties. Then $Z\subset X$ is $p^\apos$. In particular, for all locally free sheaves $\eF$ on $X$ holds: 
\begin{equation}\label{eq:isom}
\res^X_{Z}:H^t(X,\eF)\to H^t(\hat X_Z,\eF) 
\text{ is }\;\Big\{
\begin{array}{l}
\text{- an isomorphism, for }t\les p-1,\\[.75ex] 
\text{- injective, for }t=p.
\end{array}
\Big.
\end{equation}
\end{m-proposition}

\begin{m-proof}
The completion $\hat{\cal O}_{X,z}$ of the local ring at a point $z\in Z$ is isomorphic to a ring of formal power series. Consider $\xi_1,\dots,\xi_u,\zeta_1,\dots,\zeta_v\in\eO_{X,z}$, whose images in $\hat{\cal O}_{X,z}$ are independent variables, such that $\eI_{Y,z}=\lran{\bsymb{\xi}}=\lran{\xi_1,\dots,\xi_u}$ and  $\eI_{Z,z}=\lran{\bsymb{\xi},\bsymb{\zeta}}=\lran{\xi_1,\dots,\xi_u,\zeta_1,\dots,\zeta_v}$. For $l\ges a$, a direct computation yields 
$\,\eI_{Y,z}^a\cap\eI_{Z,z}^l
=\ouset{i=a}{l}{\sum}\lran{\bsymb{\xi}}^{i}\cdot\lran{\bsymb{\zeta}}^{l-i}
=\eI_{Y,z}^a\cdot\eI_{Z,z}^{l-a},$
which implies  $(\eI_{Z,z}^l+\eI_{Y,z}^{a})/{\eI_{Y,z}^{a}}\cong\eI_{Z,z}^{l}/\eI_{Y,z}^a\cdot\eI_{Z,z}^{l-a}$.
We obtain the exact sequences: 
\begin{equation}\label{eq:la}
0\to
\frac{\eI_Y^a}{\eI_Y^{a+1}}\otimes\biggl(\frac{\eI_Z}{\eI_Y}\biggr)^{l-a}
\!\to
\frac{\eI_Z^l+\eI_Y^{a+1}}{\eI_Y^{a+1}}
\to
\frac{\eI_Z^l+\eI_Y^{a}}{\eI_Y^{a}}
\to0,\quad\forall\,l\ges a+1.
\end{equation}
The left hand-side is an $\eO_Y$-module: $\eI_Z/\eI_Y=\eI_{Z\subset Y}$ is the ideal of $Z\subset Y$ and $\eI_Y^a/\eI_Y^{a+1}=\Sym^a\eN_{Y/X}^\vee$. 

Let $\eF$ be a locally free sheaf on $X$. The $p^\pos$ property implies that there are a linear function $l(k)=\cst{}_1\!\cdot\,k+\cst{}_2$ (with $\cst{}_1,\cst{}_2$ independent of $\eF$) and integers $k_\eF,l_\eF$, with the following properties: 
$$
\begin{array}{rl}
H^t(\eF\otimes\eI_{Y}^k)=0,
&
\;\forall\,t\les p,\;\forall\,k\ges k_\eF,
\\[1ex]
H^t(\eF_Y\otimes\eI_{Z\subset Y}^{l})=0,
&
\;\forall\,t\les p,\;\forall\,l\ges l_\eF,
\\[1ex] 
H^t(\eF_Y\otimes\Sym^a\eN_{Y/X}^\vee\otimes\eI_{Z\subset Y}^{l-a})=0,
&
\;\forall\,t\les p,\;\forall\,a\les k,\;\forall\,l\ges l(k).
\end{array}
$$
For the last claim, one applies the uniform $q$-ampleness property (cf. \ref{thm:unif-q}): 
\begin{itemize}
\item[--]
There is a function $\text{linear}(r)$ such that for any locally free sheaf $\eF$ with regularity ${\rm reg}(\eF_Y)\les r$ holds: 
$\;H^t(\eF_Y\otimes\eI_{Z\subset Y}^l)=0,\;\forall\,t\les p,\; l\ges{\rm linear}(r).$
\item[--] 
If $a\les k$, then ${\rm reg}(\eF_Y\otimes\Sym^a\eN_{Y/X}^\vee)\les{\rm linear}(k)$. 
\end{itemize}
Recursively for $a=1,\ldots,k$, starting with $\frac{\eI_Z^l+\eI_Y}{\eI_Y}=\eI_{Z\subset Y}^l$, the exact sequence \eqref{eq:la} yields: 
\\[.5ex] \centerline{
$
H^t\Big(
\eF\otimes\frac{\eI_Z^l+\eI_Y^{k}}{\eI_Y^{k}}
\Big)=0,\;\forall t\les p,\;\forall\,l\ges l(k). 
$
}\\[.5ex]
Now tensor $0\to\eI_Y^k\to\eI_Z^l+\eI_Y^{k}\to\frac{\eI_Z^l+\eI_Y^{k}}{\eI_Y^{k}}\to0$ by $\eF$ and deduce: 
$$
H^t\big(\eF\otimes(\eI_Y^k+\eI_Z^l)\big)=0,
\;\forall\,t\les p,\;\;\forall\,k\ges k_\eF,\;\forall\,l\ges l(k).
$$
Note that the subschemes $Z_{l,k}$ defined by $\eI_Y^k+\eI_Z^l$ are `asymmetric' thickenings of $Z$ in $X$. The sequence of ideals $\eJ_k:=\eI_Y^k+\eI_Z^{k+l(k)}$ satisfies the property \eqref{eq:ap}: indeed, 
$$
\eJ_{k'}\subset\eI_Z^k,\;\text{for}\;k'\ges k,\quad 
\eI_Z^{m'}\subset\eJ_m,\;\text{for}\;m'\ges m+l(m).
$$
Thus $Z\subset X$ is $p^\apos$. The Lemma \ref{lm:cd-apos} implies that $\cd(X\sm Z)<\dim X-p$, hence \eqref{eq:isom} holds by \cite[Theorem III.3.4(b)]{hart-as}. 
\end{m-proof}


\section{The G3 property}\label{sct:G3}

As already mentioned, in this section we shall see that partially ample subvarieties are G3 in the ambient space, not only G2 as in Section~\ref{sct:G2}. Thus we obtain extension criteria for formal meromorphic functions in a number of situations. The key to deduce the G3 property is the following result due to Speiser. 

\begin{thm-nono}{{\rm(cf. \cite[Corollary V.2.2]{hart-as})}} 
Let $X$ be a non-singular projective variety and $Y$ a closed subscheme. The statements below are equivalent: 
\begin{enumerate}[leftmargin=5ex]
\item[\rm(i)] 
$Y$ is G3 and intersects every effective divisor on $X$;
\item[\rm(ii)] 
$Y$ is G2 and $\cd(X\sm Y)\les\dim X-2$. 
\end{enumerate}
\end{thm-nono}

\begin{m-theorem}\label{thm:G3}
Let $X$ be a non-singular projective variety and $Y\subset X$ be a $(\dim Y-1)$-ample lci subvariety. Then $Y\subset X$ is G3. 
\end{m-theorem}

\begin{m-proof}
By Propositions~\ref{prop:q1}, \ref{prop:N+cd}, $Y$ is connected, $\cd(X\sm Y)\les\dim X-2$, and the normal bundle of $Y$ is $(\dim Y-1)$-ample. Corollary~\ref{cor:G2} implies that $Y\subset X$ is G2 and it remains to apply Speiser's result.
\end{m-proof}

The difficulty in proving the $1^\pos$ property of a subvariety is to control the cohomological dimension of its complement. Our strategy is to apply Speiser's theorem that is, we show the G3 property by a direct argument. In order to achieve this goal we use a variant of a result due to Chow and B\u{a}descu-Schneider (cf. \cite{chow,bad}). 

\begin{m-asmp}
Let $\cY\subset S\times X$ be a $S$-flat family of irreducible subvarieties of a non-singular variety $X$, parametrized by an irreducible quasi-projective variety $S$. We assume that the following conditions hold: 
\begin{equation}\label{eq:SYX}
\begin{minipage}{0.9\textwidth}
\begin{itemize}
\item[(i)] 
$\rho:\cY\to X$ is dominant (we say that $\cY$ is \emph{movable}); 
\item[(ii)] 
$Y_s\cap D\neq\emptyset$, for any $s\in S$ and every effective divisor $D$ on $X$. 
\end{itemize}
\end{minipage}
\end{equation}
\end{m-asmp}

\begin{m-theorem}\label{thm:badescu} 
Let the situation be as in \eqref{eq:SYX} and suppose $X$ is algebraically simply connected. If $Y_o$ is G2 in $X$, for some $o\in S$, then it is actually G3. 
\end{m-theorem}

\begin{m-proof}
The argument is analogous to \cite[Theorem 13.4(ii)]{bad}. There is a normal, irreducible, projective variety $X'$, a G3 subvariety $Y'_o\subset X'$, and a morphism 
$$f:(X',Y'_o)\to(X,Y_o)$$ 
such that $f$ is \'etale in a neighbourhood of $Y_o$ (cf. \cite[Corollary 9.20]{bad}). We are going to prove that $f$ is \'etale everywhere, so $X'=X$ since $X$ is simply connected, hence $Y_o\subset X$ is G3. 

Let $\Delta'\subset X'$ be the ramification divisor of $f$ and $\Delta:=f_*(\Delta')$; it is an effective divisor on $X$ and we want $\Delta=\emptyset$. We show that a generic deformation of $Y_o$ avoids $\Delta$, which is impossible by our hypothesis. Let us consider the diagram: 
$$
\xymatrix@R=1.75em@C=2em{
\cY'_{\bar{\si}}\ar[r]\ar[d]^-{\vphi_{\bar{\si}}}&\cY'_\si\ar[r]\ar[d]^-{\vphi_\si}
&\cY':=\cY\times_XX'\ar[r]^-{\rho'}\ar[d]^-\vphi&X'\ar[d]^-f
\\ 
\cY_{\bar{\si}}\ar[r]\ar[d]&\cY_\si\ar[r]\ar[d]&\cY\ar[r]^-\rho\ar[d]^-\pi&X
\\ 
\bar{\si}:=\Spec(\ovl{K(S)})\ar[r]&{\si}:=\Spec({K(S)})\ar[r]&S&
}
$$
The ramification divisor of $\vphi$ is $\Delta_1:=(\rho')^{-1}(\Delta')$. It is enough to prove that $\Delta_1\cap \cY'_\si=\emptyset$; in this case, $\Delta\cap Y_s=\emptyset$ for $s\in S$ generic. Moreover, it suffices to verify $\Delta_1\cap \cY'_{\bar{\si}}=\emptyset$, that is $\vphi_{\bar{\si}}$ is \'etale. 

First we observe that, since $\cY'\subset (S\times X)\times_XX'=S\times X'$ and $\rho$ is dominant, $\cY'_\si\subset X'$ is dense, thus irreducible. Second, $Y_o\times_XY'_o\cong Y_o$ is an irreducible component of $\cY'_o=\vphi^{-1}(o)$, and $\vphi$ is \'etale in its neighbourhood. But $Y_o\times_XY'_o$ is the specialization of some irreducible component $Z$ of the fibre of $\cY'_{\bar{\si}}\to\cY_{\bar{\si}}$, so $\vphi_{\bar{\si}}:Z\to\cY_{\bar{\si}}$ is \'etale. Since $\cY'_\si$ is irreducible, the components of $\cY'_{\bar{\si}}$ are conjugate under the Galois group ${\rm Gal}\big(\ovl{K(S)}/K(S)\big)$. We conclude that $\vphi_{\bar{\si}}$ is \'etale everywhere. 
\end{m-proof}

\begin{m-theorem}\label{thm:badescu2}
Let $X$ be a non-singular, algebraically simply connected variety. Consider a family of lci subvarieties $\cY$ as in \eqref{eq:SYX}, such that $\eN_{Y_o/X}^\vee$ is not pseudo-effective for some $o\in S$. Then $Y_o$ is a $1^\pos$ subvariety of $X$. \\ 
Hence, sufficiently general movable subvarieties of rationally connected varieties are $1^\pos$. 
\end{m-theorem}

\begin{m-proof}
Since the normal bundle $\eN_{Y_o/X}$ is $(\dim Y_o-1)$-ample, Corollary~\ref{cor:G2} implies that $Y_o$ is G2 in $X$; the previous theorem shows that it is actually G3. It also intersects every divisor, so Speiser's result yields $\cd(X\sm Y_o)\les\dim X-2$. We deduce that $Y_o$ is $1^\pos$, by Proposition~\ref{prop:N+cd}. 

The last claim follows from the fact that rationally connected varieties are algebraically simply connected (cf. \cite[Corollary 4.18]{debr}). 
\end{m-proof}

In section~\ref{sct:expl} we shall see that the theorem applies to quasi-homogeneous varieties (in contrast to~\cite[Section~13]{bad}). This yields many new examples of G3 subvarieties.

\subsection{Intersections with divisors}\label{ssct:int-div}

Here we discuss two fairly independent situations where the condition \eqref{eq:SYX}(ii) is satisfied:
\begin{enumerate}
\item[$\bullet$]
every effective divisor on $X$ is semi-ample that is, a multiple is globally generated;
\item[$\bullet$] 
the family $\cY$ is strongly movable.
\end{enumerate}


\subsubsection{Minimal Mori dream spaces}\label{sssct:mds}
They are the prototype of varieties whose effective divisors are semi-ample (cf. \cite[Proposition 1.11]{hu-keel}). Examples include Fano varieties, numerous toric and spherical varieties, and also GIT-quotients.  

\begin{m-corollary}\label{cor:mori-g3}
Let $X$ be a smooth, algebraically simply connected, projective variety, such that every effective divisor on it is semi-ample. Suppose $Y$ is a movable lci subvariety of $X$, such that its normal bundle is $(\dim Y-1)$-ample. Then $Y$ is G3, actually $1^\pos$, in $X$. 
\end{m-corollary}

\begin{m-proof} 
In order to apply \ref{thm:badescu2}, it is enough to show that $Y$ intersects every effective divisor $D$ on $X$ non-trivially. Suppose there is $D$ such that $Y\cdot D=0$. Since $mD$ is globally generated for some $m\ges1$, $Y$ must be contained in a fibre $F$ of $X\srel{\phi}{\to}|mD|$; the same holds for all the deformations of $Y$. Note that $\dim |mD|\ges1$. 

Since $Y$ is movable, after possibly replacing $Y$ by a deformation, we may assume that $F$ is a regular fibre of $\phi$. Then $\eN_{F/X\rst_Y}$ (which is trivial, of rank at least one) is a quotient of $\eN_{Y/X}$ (which is $(\dim Y-1)$-ample). This is impossible.
\end{m-proof}


\subsubsection{Strongly movable families}\label{sssct:strg-mov}

The concept of strongly movable subvarieties was introduced by Voisin, in the attempt to geometrically characterize big subvarieties of projective varieties (cf. \cite[Section~2]{voisin-coniv}). 

\begin{m-notation}\label{not:syx}
Let $\cY\srel{(\pi,\rho)}{\subset} S\times X$ be a family of lci subvarieties of $X$, with $\rho$ dominant; then $\rho(\cY)$ contains an open subset $O$ of $X$. 
The incidence variety $\Si$ is the irreducible component of $(\pi,\pi)(\cY\times_X\cY)$ which contains the diagonal; $\pi$ is a projective morphism, so $\Si\subset S\times S$ is closed. One obtains the diagram:
\begin{equation}\label{eq:rho}
\xymatrix@C=5em{
\cY_\Si:=\Si\times_{S\times S}(\cY\times_X\cY)\ar[r]\ar@/^3ex/[rr]^-{\rho_\Si}\ar[d]&
\cY\times_X\cY\ar[r]_-{(\pi,\rho)}\ar[d]^-{(\pi,\pi)}&S\times X.
\\ 
\Si:={\rm Image}(\pi,\pi)\ar[r]^-{\iota}&S\times S&
}
\end{equation}
For $o\in S$, we denote $\Si_o:=\iota^{-1}(o,S)$ and $\rho_o:=\rho_{\Si\rst_{\Si_o}}$. 
\end{m-notation}

\begin{m-definition}\label{def:mov+conn}
\begin{enumerate}
\item[(i)]
The family $\cY$ is \emph{strongly movable}, if $\rho_\Si$ is dominant.

\item[(ii)] 
The variety $X$ is $\cY$-\emph{chain connected in codimension one}, if there is an open subset $O\subset\rho(\cY)$ satisfying the following properties:
$$
\begin{array}{l}
\bullet\;\codim_X(X\sm O)\ges2,\;\text{and}\\[1.5ex] 
\bullet\;\forall x,x'\in O,\;\exists s_0,\dots,s_n\in S,\,
x\in Y_{s_0},\,x'\in Y_{s_n},\;Y_{s_{j-1}}\cap Y_{s_j}\neq\emptyset,\;j=1,\dots,n.
\end{array}
$$
For simplicity, we call such a sequence $(Y_{s_0},\dots,Y_{s_n})$ a \emph{$\cY$-chain}. 
\end{enumerate}
\end{m-definition}
Obviously, if $\cY$ is strongly movable, then there is an open subset of $X$ whose points are connected by $\cY$-chains. We require that its complement in $X$ has codimension at least two.

\begin{m-lemma}\label{lm:1+}
Let the notation be as in \ref{not:syx} and suppose $\cY$ is strongly movable. Then, for generic $o\in S$, the normal bundle $\eN_{Y_o/X}$ is $(\dim Y_o-1)$-ample. 
\end{m-lemma}

\begin{m-proof}
We must find a movable morphism $C\srel{\vphi}{\to}Y_o$, an ample line bundle $\eL_C\in\Pic(C)$, and a movable homomorphism $\eL_C\to\vphi^*\eN_{Y_o/X}$. A tangent vector $\xi\in\eT_{S,o}$ induces an infinitesimal deformation $\hat v_\xi\in H^0(Y_o,\eN_{Y_o/X})$ that is, a homomorphism 
$$
\hat v_\xi:\eO_{Y_o}\to\eN_{Y_o/X}.
$$
Henceforth we restrict our attention to $\xi\in\eT_{\Si_{o},o}\subset\eT_{S,o}$. 
\smallskip

\nit\underbar{\textit{Claim~1}}\quad 
The vanishing locus of $\hat v_\xi$ is a non-empty, proper subset of $Y_o$. Moreover, for $\xi\in\eT_{\Si_o,o}$ variable, the vanishing loci of $\hat v_\xi$ cover an open subset of $Y_o$. \\ 
The vector $\xi$ is determined by an arc $\Spec\big(\kk\ldbrack\eps\rdbrack\big)\srel{h}{\to}\Si_o\,$ through $o$ on $\Si_o$. The defining property of $\Si_o$ implies that there is $y_\eps\in Y_o\cap Y_{h(\eps)}$, $\disp\lim_{\veps\to0} y_\eps=y\in Y$. By differentiation, we deduce that the infinitesimal deformation $\hat v_\xi$ vanishes at $y$, because $Y_o$ is deformed in a tangential direction at that point. 

We claim that $\hat v_\xi\neq0$, for generic $\xi$, and the vanishing loci cover an open subset of $Y_o$. Both statements follow from the strongly movability of $\cY$, which implies that $\cY_{\Si_o}\srel{\rho_o}{\to}X$ is dominant (recall that $o\in S$ is generic). Indeed, on one hand, the differential of $\rho_o$ at a generic point $y\in Y_o$ is surjective that is, $\eT_{\Si_o,o}\srel{\rd\rho_{o,y}}{\lar}\eN_{Y_o/X,y}$, $\,\xi\mt\hat v_{\xi,y}$, is surjective. On the other hand, $\rho_o^{-1}(Y_o)\to Y_o$ is dominant too, so the points $y$ as above cover an open subset of $Y_o$. \smallskip

\nit\underbar{\textit{Claim~2}}\quad 
Now let $C\subset Y_o$ be a complete intersection curve which intersect the zero locus of $\hat v_\xi$ properly. By Claim~1, such curves are movable. Moreover, $\hat v_\xi$ extends to a pointwise injective homomorphism $\eL_C\subset\eN_{Y_o/X\rst_C}$, where $\eL_C$ is an ample line bundle.  The latter is movable too, because $\rd\rho_{o,y}$ is surjective at the generic point $y\in Y_o$.
\end{m-proof}

\begin{m-lemma}\label{lm:conn}
Let the notation be as in \ref{not:syx} and suppose $X$ is $\cY$-chain connected in codimension one. Then, for any $s\in S$, $Y_s$ intersects every effective divisor $D$ on $X$ and the intersection $D\cdot Y_s$ is numerically non-trivial. 
\end{m-lemma}

\begin{m-proof}
We fix a generic point $o\in S$ and assume that there is an irreducible effective divisor $D$ in $X$, such that $D\cap Y_o=\emptyset$. Note that $D\cap O\neq\emptyset$, so $\cal D:=\overline{\rho^{-1}(D\cap O)}\subset\rho^{-1}(D)$ is a divisor on $\cY$. \smallskip

\nit\underline{\textit{Claim~1}}\quad $D\cap Y_s=\emptyset$ or $Y_s\subset D$, for all $s\in S$.\\ 
Otherwise, for some $s\in S$, the intersection $D\cap Y_s$ is proper, thus $D\cdot Y_s$ is numerically non-trivial. (Its pairing with a general intersection of complementary dimension is non-empty.) The deformation invariance of the intersection product yields $D\cap Y_o\neq0$, a contradiction. 

Note that, since $\pi$ is projective, $\pi(\cal D)$ is closed in $S$. The claim implies:
$$
\left\{\begin{array}{l}
\cal D\cap\rho^{-1}(Y_s)=\emptyset\\[.5ex]\text{or}\\[.5ex]\rho^{-1}(Y_s)\subset\cal D
\end{array}\right.
\quad\Rightarrow\quad
\left\{\begin{array}{l}
\pi(\cal D)\cap\pi(\rho^{-1}(Y_s))=\emptyset\\[.5ex]\text{or}\\[.5ex]\pi(\rho^{-1}(Y_s))\subset\pi(\cal D).
\end{array}\right.
$$
For the last implication, observe the following: 
$$
t\in\pi(\cal D)\cap\pi(\rho^{-1}(Y_s))
\quad\Rightarrow\quad 
D\cap Y_t\neq\emptyset,\;\;\text{so}\;Y_t\subset D
\quad
\srel{\kern-3ptY_s\cap Y_t\,\subset\,Y_t}%
{=\kern-.5ex=\kern-.5ex=\kern-.5ex=\kern-.5ex=\kern-.5ex\Longrightarrow}\quad 
Y_s\subset D.
$$ 

\nit\underline{\textit{Claim~2}}\quad $\pi(\cal D)$ contains an open subset of $S$, so $\pi$ is surjective.\\ 
We use the chain-connectedness of $\cY$: since $D\cap O\neq\emptyset$, there are $o=s_0,s_1,\dots,s_n\in S,$  
$$
Y_{s_{j-1}}\cap Y_{s_j}\neq\emptyset\;\;\textrm{and}\;\;D\cap Y_{s_n}\neq\emptyset.
$$
The previous claim implies $Y_{s_n}\subset D$, so $D\cap Y_{s_{n-1}}\neq\emptyset$, \textit{etc}. Inductively, we deduce $Y_o\subset D$, so the generic point $o\in S$ belongs to $\pi(\cal D)$. Then $\pi(\cal D)=S$, so it holds $\cal D=\pi^{-1}(\pi(\cal D))=\cY$, a contradiction.   

Finally, the intersection $D\cap Y_o$ is proper, hence numerically non-trivial, because otherwise $Y_o\subset D$, for generic $o\in S$, which is impossible.   
\end{m-proof}

\begin{m-theorem}\label{thm:strg-mov}
Let $\cY$ be a family of lci subvarieties of the smooth projective variety $X$. We assume that the following properties are satisfied: 
\begin{enumerate}[leftmargin=5ex]
\item[\rm(a)]
$X$ is algebraically simply connected (e.g. $X$ is rationally connected);
\item[\rm(b)]
$\cY$ is strongly movable;
\item[\rm(c)]
$X$ is $\cY$-chain connected in codimension one.
\end{enumerate}
Then any member $Y$ of $\cY$ is G3 in $X$. 
\end{m-theorem}

\begin{m-proof}
We apply~\ref{lm:1+} and \ref{lm:conn}: $\eN_{Y/X}$ is $(\dim Y-1)$-ample, hence $Y$ is G2 in $X$. Moreover, $Y$ intersects every divisor, so $Y\subset X$ is G3. 
\end{m-proof}

\begin{m-remark}
The smoothness of $X$ can be weakened: it is enough if the generic member of the family $\cY$ is contained in the smooth locus of $X$.
\end{m-remark}


\section{A connectedness problem}\label{sct:connect}

\begin{m-notation}\label{fh}
Let $f:V\to X$ be a morphism between irreducible projective varieties and $Y\subset X$ be a lci subvariety, with $\codim_XY<\dim f(V)$. We denote 
$$
q:=\dim f(V)+\dim Y-\dim X-1=\dim f(V)-\codim_XY-1,
$$
and assume that $q\ges0$.
\end{m-notation}

The goal of this section is to apply the ideas developed so far to the following:

\begin{conj-nono}{} 
If the normal bundle $\eN_{Y/X}$ is ample, then $f^{-1}(Y)$ is connected. 
\end{conj-nono}

The conjecture is due to Fulton-Hansen (cf. \cite[pp. 161]{fult+hans}). Hartshorne raised a similar problem for $f$ a closed embedding, concerning the non-emptiness of $Y\cap V$ (cf. \cite[Ch.~III, Conjecture~4.5]{hart-as}). As stated, the conjecture is not true: there are counterexamples constructed by Hartshorne, see Remark~\ref{rmk:connect}(ii) below.

The conjecture holds for products of projective spaces, flag varieties, low codimensional subvarieties of projective spaces, by work of Fulton-Hansen, Debarre, B\u{a}descu \cite{fult+hans,hans,debr2,bad-debarre,bad-BL}. Faltings proved the statement for homogeneous varieties (cf. \cite[Korollar, pp. 148]{falt-homog}). The common feature of these approaches is to focus on the properties of the diagonal $\Delta_X\subset X\times X$. Moreover, the methods are strongly adapted to homogeneous varieties. The problem was also studied in \cite{ran}, where the emphasis is rather on the movability of $Y$ in $X$. 

The approach inhere is based on two ingredients: Theorem~\ref{thm:G3} and the following result due to Hironaka-Matsumura: 
\begin{thm-nono}{{\rm(cf. \cite[Theorem~2.7]{hir+mats})}}
Suppose $V\srel{f}{\to}X$ is \emph{surjective} and $Y\subset X$ is G3. Then $f^{-1}(Y)\subset V$ is G3, hence connected.
 \end{thm-nono}
Hartshorne's example shows that the G3-condition is optimal. This motivates our strategy: we impose positivity conditions on $Y$ and analyse how they are inherited by the pre-image. 


\begin{m-theorem}\label{thm:fh}
Let the notation be as in~\ref{fh} and assume that $Y\subset X$ is $q$-ample. 
\begin{enumerate}[leftmargin=5ex]
\item[\rm(i)] Suppose $f$ is an embedding. Then the following statements hold:
\begin{enumerate}
\item[\rm(a)]
If $V$ is Cohen-Macaulay, then $Y\cap V$ is non-empty and connected. 
\item[\rm(b)]
If $V$ is smooth and $Y\cap V\subset V$ is lci, then $Y\cap V$ is G3 in $V$.  
\end{enumerate}
\item[\rm(ii)] Let $f:V\to X$ be a morphism.  
\begin{enumerate}
\item[\rm(a)]
If the Stein factorization of $f$ is Cohen-Macaulay, then $f^{-1}(Y)$ is non-empty and connected. 
\item[\rm(b)]
The pre-image $f^{-1}(Y)$ is G3 in $V$ if one of the following conditions is satisfied:
\begin{enumerate}
\item[$\circ$]
$V$ is smooth, $f$ is flat, and $Y\cap f(V)$ is lci in $f(V)$;
\item[$\diamond$]
$f(V)$ is smooth and $Y\cap f(V)$ is lci in $f(V)$.
\end{enumerate}
\end{enumerate}
\end{enumerate}
\end{m-theorem}

\begin{m-proof}
(i)(a) We claim that $Y\cap V$ is non-empty and $\big(\dim(Y\cap V)-1\big)$-ample in $V$: indeed, 
$\;\cd(V\sm Y\cap V)\les\cd(X\sm Y)\les\codim_XY+q-1=\dim V-2$ (cf. Proposition~\ref{prop:EY}). 

Thus $Y$ intersects $V$, actually $\dim(Y\cap V)\ges1$; otherwise $\cd(V\sm Y\cap V)=\dim V-1$, a contradiction. The universality property of the blow-up yields the Cartesian diagram: 
$$
\xymatrix@R=1.5em{
\Bl_{Y\cap V}(V)\ar@{^(->}[r]\ar[d]&\Bl_{Y}(X)\ar[d]
\\ 
V\ar@{^(->}[r]&X.
}
$$
Since the exceptional divisor $E_{Y}$ is $(\dim V-2)$-ample, the same holds for $E_{Y\cap V}$. Proposition~\ref{prop:connect} implies that $Y\cap V$ is connected. 

\nit(b) It's a direct application of Theorem~\ref{thm:G3}.\medskip

\nit(ii) Let $Z:=f(V)$ be the reduced image and $\bar V:=\Spec(f_*\eO_V)$ the Stein factorization of $f$.\\ 
\nit(a) As before, it holds $\cd(Z\sm Y\cap Z)\les\cd(X\sm Y)\les\dim Z-2$. Since $\bar V\srel{\bar f}{\to}Z$ is finite and $\bar V$ is Cohen-Macaulay, it follows $\cd(\bar V\sm{\bar f}^{-1}(Y))\les\dim\bar V-2$, thus ${\bar f}^{-1}(Y)$ is connected. But $V\to\bar V$ has connected fibres, so $f^{-1}(Y)$ is connected.

\nit(b$\,\circ$) In this case, $f^{-1}(Y)\subset V$ is lci and $f$ is equidimensional. The induced morphism $\tld f:\Bl_{f^{-1}(Y)}(V)\to\Bl_{Y\cap Z}(Z)$ is still equidimensional; since $E_{Y\cap Z}$ is $(\dim Z-2)$-ample, $E_{f^{-1}(Y)}$ is $(\dim V-2)$-ample (cf. proof of ~\ref{prop:pull-back}). The claim follows from Theorem~\ref{thm:G3}. 

\nit(b$\,\diamond$) By (i)(b), $Y\cap Z\subset Z$ is G3, so $f^{-1}(Y)\subset V$ is G3.
\end{m-proof}

\begin{m-remark}\label{rmk:connect}
(i) If $Y$ is an ample subvariety of $X$, the condition $q\ges0$ in~\ref{fh} becomes the usual one: $\;\dim f(V)>\codim_XY$. Note that, in this case, $\eN_{Y/X}$ is ample (cf.~\ref{prop:N+cd}). 

\nit(ii) We pointed out that, as stated, the conjecture is false: \cite[pp. 199]{hart-as} is the example of an \'etale morphism $V\to X$ and positive dimensional, smooth subvariety $Y\subset X$ with ample normal bundle and disconnected pre-image. The condition (iib$\,\circ$) is satisfied and, moreover, $Y$ is G2 but not G3 in $X$. What (necessarily) fails is the partial ampleness of $Y$. 
\end{m-remark}

Hence the positivity conditions on $Y$ are necessary. However, one can weaken the transversality and regularity assumptions through deformations. The method is suited for smoothable triples $(Y,V,f)$, which can be moved into general relative position.

\begin{m-lemma}{\rm(cf. \cite[Sect.~2]{fult+hans})}\label{lm:deform}
Let $S$ be a normal, irreducible, quasi-projective variety,
\begin{enumerate}[leftmargin=5ex]
\item[$\bullet$]
$\cY=\{Y_s\}_{s\in S}\srel{(\pi,\iota)}{\subset}S\times X$ a flat family of $c$-codimensional subvarieties,
\item[$\bullet$]
$\cal V=\{V_s\}_{s\in S}$ a $S$-flat family of projective varieties, 
$\cal V\srel{f=(f_s)_{s\in S}}{-\kern-1ex-\kern-1ex\lar}\!X$ a morphism.
\end{enumerate}
Consider $\cal W=(W_s)_{s\in S}:=\cY\times_{S\times X}\cal V\subset\cal V$; assume that generically $W_s\subset V_s$ is connected and $c$-co\-dimensional, for $s\in\cal U\subset S$. 
For $o\in S$, suppose moreover that it holds:
$$
\forall\,y_o\in Y_o\cap f(V_o),\; \exists\,\Spec(\kk\ldbrack\eps\rdbrack)\srel{u}{\to}\cY\times_{S\times X}f(\cal V)\subset\cY,
\;u(0)=y_o,\;u(\eps)\in\pi^{-1}(\cal U).
$$
Then $f_o^{-1}(Y_o)\subset V_o$ is connected.
\end{m-lemma}


\section{Examples of \textit{q}-ample subvarieties}\label{sct:expl}

Partially ample subvarieties are ubiquitous. We shall discuss several classes of examples:
\begin{enumerate}
\item subvarieties of almost homogeneous varieties;
\item zero loci of sections in globally generated vector bundles;
\item sources of Bialynicki-Birula decompositions, corresponding to actions of the multiplicative group. 
\end{enumerate}

\subsection{Subvarieties of (almost) homogeneous varieties}\label{ssct:subvar-homog}

Suppose $G$ is a connected algebraic group with identity element $e$ and consider the morphism: 
$$
\gamma:G\times G\to G,\quad\gamma(g',g):=g'g^{-1}.
$$ 
Let $X$ be a smooth $G$-variety with an open orbit $O$. So, if $G$ is linear, then $X$ is automatically rational. We denote by $\mu:G\times X\to X$ the action and, for $x\in X$, let $\mu_x(\cdot):=\mu(\cdot,x)$. 

\begin{m-definition}\label{def:gener} 
Suppose $Y\subset X$ is an irreducible subvariety which intersects $O$. We consider the following objects: 
\begin{itemize}[leftmargin=5ex]
\item[$\bullet$] 
For $y_o\in Y\cap O$, let $G_{Y,y_o}:=\mu^{-1}_{y_o}(Y)$; it is a closed subvariety of $G$. Denote $G_Y$ the subgroup generated by $G_{Y,y_o}$; it is closed in $G$, independent of $y_o\in Y\cap O$.
\item[$\bullet$] 
Define $S_Y:=\gamma(G_{Y,y_o},G_{Y,y_o})$; it is constructible in $G$, contains $G_{Y,y_o}$, and independent of $y_o\in Y\cap O$. In fact, $S_Y$ consists of those elements of $G$ which send some point of $Y\cap O$ to another point of $Y\cap O$. 
\end{itemize}
We say that $Y$ \emph{generates} $X$ if $G_{Y}=G$. We say that $Y$ \emph{strongly generates} $X$ if $\mu(S_Y,Y\cap O)$ is open in $O$. 
\end{m-definition}

Note that, for homogeneous $X=G/P$, the subset $S_Y$ is closed in $G$. Indeed, if $y_o=e$, then $P\subset G_{Y,y_o}$ and $\mu$ factorizes through $(G\times G)/P$, where $P$ acts diagonally. But the morphism $(G\times G)/P\to G$ is projective, so the image of $(G_{Y,y_o}\times G_{Y,y_o})/P$ is closed in $G$.

\begin{m-lemma}\label{lm:strg-gen-gen}
If $Y$ strongly generates $X$, then it generates $X$.
\end{m-lemma} 

\begin{m-proof}
Indeed, the stabilizer $H\subset G$ of $y_o$ is contained in $G_{Y,y_o}$ and 
$$
\mu(S_Y,Y\cap O)=\mu\big((G_{Y,y_0}\cdot G_{Y,y_o}^{-1}\cdot G_{Y,y_o}),y_o\big)
\cong 
(G_{Y,y_0}\cdot G_{Y,y_o}^{-1}\cdot G_{Y,y_o})/H.
$$ 
If $\mu(S_Y,Y\cap O)$ is open in $O\cong G/H$, then $G_{Y,y_0}\cdot G_{Y,y_o}^{-1}\cdot G_{Y,y_o}$ is $H$-invariant and open in $G$, so $G_Y=G$. 
\end{m-proof}

We remark that, by taking $G$-translates of $Y$, we retrieve the situation studied in \ref{sssct:strg-mov}:
\begin{equation}\label{eq:strg-gen}
\xymatrix{
\cY:=G\times Y\ar[r]^-\mu\ar[d]_-\pi&X\\ 
G&
}
\end{equation}

\begin{m-lemma}\label{lm:strg-gen}
Suppose $Y$ strongly generates the almost homogeneous variety $X$. Then the following properties hold:
\begin{enumerate}[leftmargin=5ex]
\item[\rm(i)]
the family \eqref{eq:strg-gen} is strongly movable;
\item[\rm(ii)]
If $\codim_X(X\sm O)\ges2$, then $X$ is $\cY$-chain connected in codimension one.
\end{enumerate}
\end{m-lemma}

\begin{m-proof}
(i) Note that $S_Y$ coincides with $\Si_e$ defined in \ref{not:syx} (for $o=e\in G$). Moreover, for $g\in G$, $S_{gY}=gS_Yg^{-1}$ and it holds: $\;\mu(S_{gY},gY\cap O)=g\mu(S_Y,Y\cap O).$
By hypothesis, the right hand-side is open in $X$, hence $\rho_\Si$ in \eqref{eq:rho} is dominant.

\nit(ii) Lemma~\ref{lm:strg-gen-gen} implies that $G$ is generated by $G_{Y,y_o}$. Then, for any $g\in G$, there is a sequence $e=g_0,g_1,\dots,g_n\in G_{Y,y_o}$ whose product equals $g$. It follows that 
$$
Y_0:=Y,Y_1:=g_1Y,\dots,Y_n:=g_1\dots g_nY
$$
is a $\cY$-chain connecting $y_o$ to $gy_o$. 
\end{m-proof}

\begin{m-theorem}\label{thm:G3-a-homog} 
Let $(X,O)$ be an almost homogeneous variety for the action of a linear algebraic group $G$. Suppose $Y\subset X$ is a smooth subvariety with the following properties: 
\begin{enumerate}[leftmargin=5ex]
\item[\rm(a)] 
It strongly generates $X$.
\item[\rm(b)] 
The intersection with every divisor is numerically non-trivial.\\ 
(If $\codim_X(X\sm O)\ges 2$, it suffices $Y\cap O\neq\emptyset$.)
\end{enumerate} 
Then $Y$ is $1^\pos$, in particular it is G3 in $X$.  
\end{m-theorem}

\begin{m-proof}
Since $G$ is linear, $X$ is a rational variety. The previous lemma implies that Theorem~\ref{thm:strg-mov} applies to the family $\cY=G\times Y$. 
\end{m-proof}

\begin{m-corollary}\label{cor:DX}
Suppose $(X,O)$ is almost homogeneous for the action of the linear algebraic group $G$. Suppose that the following conditions are satisfied: 
\begin{enumerate}[leftmargin=5ex]
\item[\rm(a)]
$\codim_X(X\sm O)\ges2$;
\item[\rm(b)]
The stabilizer of some (any) point $x_0\in O$ contains a Cartan subgroup of $G$. 
\end{enumerate}
Then the diagonal is $1^\pos$, thus G3, in the product $X\times X$. 
\end{m-corollary}
This generalizes to almost homogeneous spaces the analogous result for rational homogeneous varieties (cf. \cite[Theorem 2]{bad+schn}).

\begin{m-proof}
We apply the previous theorem to the diagonal $\Delta_X\subset X\times X$. The group $G\times G$ acts on $X\times X$ with open orbit $O\times O$, whose complement has codimension at least two, and $\Delta_X\cap(O\times O)\cong O$. It remains to verify that $\Delta_X$ is strongly generating that is, $S_{\Delta_X}\cdot O$ is open in $O\times O$. Let $H$ be the stabilizer of the point $x_o\in O$. A direct computation yields:  
$$
\begin{array}{rl}
S_{\Delta_X}&
=\{ (g,g')\in G\times G\mid\exists x_1,x_2\in O,\;(g,g')\cdot(x_1,x_1)=(x_2,x_2) \}
\\[1.5ex]&
=\{ (g,g)\cdot(e,\Ad_a(h))\mid a,g\in G,h\in H\},
\;\;(\text{$e\in G$ is the identity}). 
\end{array}
$$
Since $H$ contains a Cartan subgroup of $G$, $\uset{a\in G}{\bigcup}\Ad_a(H)$ contains an open subset of $G$, so $S_{\Delta_X}$ contains an open subset of $G\times G$. 
\end{m-proof}
Note that, instead of (b), it suffices $H\,{\cdot}\uset{a\in G}{\bigcup}\Ad_a(H)$ to contain an open subset of $G$.


\subsubsection*{The case of homogeneous varieties}\label{sssct:IP}

In the previous paragraph we obtained the $1^\pos$ property for strongly generating subvarieties. However, one would like to have examples of subvarieties with stronger positivity properties. The theorem below is simply a reformulation of results due to Faltings, Barth-Larsen. 

\begin{m-theorem}\label{thm:falt-homog}
Let $X=G/P$ be a rational homogeneous variety, with $G$ a semi-simple linear group. Denote by $\ell$ the minimal rank of the simple factors of $G$. Let $Y\subset X$ be a smooth subvariety of codimension $\delta$. The following statements hold: 
\begin{enumerate}[leftmargin=5ex]
\item[\rm(i)]{\rm(cf. \cite[Satz~5, Satz~7]{falt-homog})} $Y$ is $(\ell-2\delta+1)^\pos$. \\ 
Hence, for $X=\mbb P^n$, $Y$ is $(2\dim Y-n+1)^\pos$ 
(cf. \cite{barth}, \cite[Theorem~3.2.1]{laz}). 
\item[\rm(ii)]{\rm(}cf. \cite[Theorem 7.1]{ottm}{\rm)}
Let $Y\subset\mbb P^n$ be a smooth subvariety. Then holds: 
$$
Y\mbox{ is }p^\pos\;\Leftrightarrow\;
\res_Y^t:H^t(\mbb P^n;\mbb Q)\to H^t(Y;\mbb Q) 
\mbox{ is an isomorphism},\;\forall\,t<p. 
$$
\end{enumerate}
\end{m-theorem}

\begin{m-proof}
(i) Indeed, we have $\cd(X\sm Y)\les\dim X-\ell+2\delta-2$ and $\eT_X$ is $(\dim X-\ell)$-ample. Since $\eN_{Y/X}$ is a quotient of $\eT_X$, we conclude by Proposition~\ref{prop:N+cd}. 

\nit(ii) In this case, $\eN_{Y/\mbb P^n}$ is ample. By \cite[Theorem 4.4, 2.13]{ogus}, the inequality $\cd(\mbb P^n\sm Y)<n-p$ holds if and only if $\res_Y^t$ is an isomorphism and the local cohomological dimension of $Y\subset\mbb P^n$ is at most $n-p$. Since $Y$ is smooth, the latter equals $\codim_{\mbb P^n}Y=n-\dim Y$, so the latter condition is satisfied. 
\end{m-proof}


\subsection{Zero loci of sections in globally generated vector bundles}\label{ssct:glob-gen}

Throughout this section, $\eN$ is a vector bundle of rank $\nu$ on the smooth projective variety $X$.

\subsubsection{$q$-ample vector bundles}\label{sssct:q-vb}

A simple method to produce $q$-ample, lci subvarieties is by taking zero loci of $q$-ample vector bundles. 

\begin{m-proposition} \label{prop:q21} 
Suppose $\eN$ is $q$-ample and $Y$ is the zero locus of a \emph{regular} section in it. Then $Y\subset X$ is a $q$-ample subvariety; in particular, if $q<\dim X-\nu$, then $Y$ is G3 in $X$. 
\end{m-proposition}

\begin{proof} 
We verify the condition \eqref{eq:q12} for the zero locus $Y$ of $s\in H^0(X,\eN)$ and the (arbitrary) vector bundle $\eF$ on $X$. Since $s$ is regular, $Y$ is lci in $X$, $\codim_X(Y)=\nu$, so $\Bl_Y(X)$ is Gorenstein, and we have the resolution (cf. \cite[Theorem 3.1]{bu+ei})
\begin{equation}\label{eq:koszul-m}
0\to L^\nu_m(\eN^\vee)\to\ldots\to L^j_m(\eN^\vee)\to\ldots\to 
\Sym^m(\eN^\vee)\srel{s^m\ort}{-\kern-1ex-\kern-1ex\lar}\eI_Y^m\to 0,\;\;\forall\,m\ges 1.
\end{equation}
The vector bundles $L^j_m(\eN^\vee)$, $1\les j\les\nu$, are defined as follows:
\\ \centerline{$
L^j_m(\eN^\vee):=
\Img\Bigl(
\Sym^{m-1}(\eN^\vee)\otimes\overset{j}{\hbox{$\bigwedge$}}\,\eN^\vee 
\srel{\phi^j_m}{-\kern-1ex\lar} 
\Sym^{m}(\eN^\vee)\otimes\overset{j-1}{\hbox{$\bigwedge$}}\,\eN^\vee 
\Bigr).
$}
The general linear group is linearly reductive and $\phi^j_m$ is equivariant, so $L^j_m(\eN^\vee)$ is actually a direct summand in $\Sym^m(\eN^\vee)\otimes\overset{j-1}{\hbox{$\bigwedge$}}\,\eN^\vee$. 
Since $\eN$ is $q$-ample, one has 
\\ \centerline{$
H^{t+j-1}(X,\eF\otimes\oset{j-1}{\bigwedge}\eN^\vee\otimes\Sym^m\eN^\vee)=0,\quad j=1,\dots,\nu,
$}\\[.5ex]
for $\,t+\nu-1\les\dim X-q-1$ and $m\gg0$. We deduce: 
$$
H^t(X,\eF\otimes\eI_Y^m)=0,\quad \forall\,0\les t\les\dim Y-q,\;m\gg0, 
$$
hence \eqref{eq:q12} is satisfied.
\end{proof}

We are going to see that the criterion is not optimal (cf. \ref{rmk:sommese-weak}).


\subsubsection{Globally generated vector bundles}\label{ssct:x-b}

In the remaining part of the section, we assume that $\eN$ is \emph{globally generated}. Suppose $Y\subset X$ is lci of codimension $\delta$, and the zero locus of a section $s\in\Gamma(\eN):=H^0(X,\eN)$. We \emph{do not require} $s$ to be regular, so we allow $\delta<\nu$. 

In this context, the situation \ref{prop:x-b} arises as follows: the blow-up fits into 
\begin{equation}\label{eq:tld-x}
\xymatrix@R=1.5em@C=1.5em{
\tld X\ar@{^(->}[r]\ar[d]_-\pi\ar@<-5pt>[rrd]_-\phi
&
\mbb P(\eN)
{=}\,
\mbb P\Bigl(
\mbox{$\overset{\nu-1}{\bigwedge}$}\eN^\vee\otimes\det(\eN)
\Bigr)\ar@{^(->}[r]
&
X\times\mbb P\Bigl(
\mbox{$\overset{\nu-1}{\bigwedge}$}\Gamma(\eN)^\vee
\Bigr)\ar[d]
\\ 
X&&\mbb P:=\mbb P\Bigl(
\mbox{$\overset{\nu-1}{\bigwedge}$}\Gamma(\eN)^\vee
\Bigr),
}
\end{equation}
and it holds 
\begin{equation}\label{eq:o1}
\eO_{\tld X}(E_Y)=\eO_{\mbb P(\eN)}(-1)\big|_{\tld X}=
\bigl(\det(\eN)\boxtimes\eO_{\mbb P}(-1)\bigr)\big|_{\tld X}.
\end{equation}

\begin{m-proposition}\label{prop:p}
Suppose $\det(\eN)$ is ample. If the dimension of the generic fibre of $\phi$ (over its image) is $p+1$, then $\eO_{\tld X}(E_Y)$ is $\dim\phi(\tld X)$-positive, and $Y$ is $p^\pos$. 
\end{m-proposition}

\begin{proof}
The assumptions of \ref{prop:x-b} are satisfied. 
\end{proof}



\subsubsection{Special subvarieties of the Grassmannian}\label{sssct:spec-grass}

Let $W\subseteq\Gamma(\eN)$ be a vector subspace which generates $\eN$, $\dim W=\nu+u+1$. It is equivalent to a morphism $f:X\to\Grs(W;\nu)$ to the Grassmannian of $\nu$-dimensional quotients of $W$; then $\det(\eN)$ is ample if and only if $\vphi$ is finite onto its image. 
Henceforth we restrict our attention to $X=\Grs(W;\nu)$; it is naturally isomorphic to the variety $\Grs(u+1;W)$ of $(u+1)$-dimensional subspaces of $W$. 

Let $\eN$ be the universal quotient bundle. The morphism $\phi$ in \eqref{eq:tld-x} is explicit: 
\begin{equation}\label{eq:q}
\mbb P(\eN)\to\mbb P,\quad 
(x,\lran{e_x})\mt 
\det(\eN_x/\lran{e_x})^\vee\subset
\oset{\nu-1}{\bigwedge}\eN_x^\vee\subset
\oset{\nu-1}{\bigwedge}W^\vee.
\end{equation}
($\lran{e_x}$ stands for the line generated by $e_x\in\eN_x$, $x\in\Grs(W;\nu)$.) 
The restriction to the Grassmannian corresponds to the commutative diagram  
\begin{equation}\label{eq:wn}
\xymatrix@R=1.75em{
0\ar[r]&\eO_{\Grs(W;\nu)}\ar[r]^-{s}\ar@{=}[d]&
W\otimes\eO_{\Grs(W;\nu)}\ar[r]\ar@{->>}[d]^-{\;\beta}&
W/\lran{s}\otimes\eO_{\Grs(W;\nu)}\ar[r]\ar@{->>}[d]&0
\\ 
&\eO_{\Grs(W;\nu)}\ar[r]^-{\beta s}&\eN\ar[r]&\eN/\lran{\beta s}\ar[r]&0.
}
\end{equation}
Thus $\phi$ is the desingularization of the rational map 
\begin{equation}\label{eq:q-grs}
g_s:\Grs(W;\nu)\dashto\Grs(W/\lran{s};\nu-1),\quad 
[W\surj N]\mt [W/\lran{s}\;\surj N/\lran{\beta s}],
\end{equation}
followed by the Pl\"ucker embedding of $\Grs(W/\lran{s};\nu-1)$. The indeterminacy locus of $\phi$ is $\Grs(W/\lran{s};\nu)$; hence the latter is $u^\pos$ in $\Grs(W;\nu)$.

The discussion generalizes to special Schubert subvarieties. For $\ell\les\nu$, fix an $\ell$-dimensional subspace $\Lambda_\ell\subset W$ and define:  
$$
Y_\ell:=\{U\in\Grs(u+1;W)\mid U\cap\Lambda_\ell\neq0\}.
$$
The cycles $Y_\ell$, $\ell=1,\dots,\nu$, generate the Chow (cohomology) ring of $X$. In fact, $Y_\ell$ has codimension $\nu-\ell+1$ and it is the vanishing locus of the section: 
$$
s_\ell:\eO\cong\det(\Lambda_\ell\otimes\eO)\to\oset{\ell}{\bigwedge} W\otimes\eO\to\oset{\ell}{\bigwedge}\eN.
$$ 
We are in the situation described in \ref{prop:x-b}. The diagram \eqref{eq:tld-x} corresponds to the rational map 
$$
\phi:\Grs(u+1;W)\dashto\Grs(u+1;W/\Lambda_\ell),\quad U\mt(U+\Lambda_\ell)/\Lambda_\ell, 
$$
followed by a large Pl\"ucker embedding; its indeterminacy locus is precisely $Y_\ell$. Since $\phi$ is surjective, we deduce that $Y_\ell\subset X$ is $\big(\ell(u+1)-1\big)^\pos$.

\begin{m-remark}\label{rmk:sommese-weak}
\nit(i) Propositions~ \ref{prop:p} and \ref{prop:q21} (see also \ref{prop:q2}) deal with complementary situations: $\eO_{\tld X}(E_Y)$ is relatively ample for some morphism, while $\eO_{\mbb P(\eN^\vee)}(1)$ is the pull-back of an ample line bundle.\smallskip

\nit(ii) 
The criterion \ref{prop:q21} is not optimal: Proposition~\ref{prop:q2} implies that, for $X=\Grs(\nu+u+1;\nu)$, the universal quotient bundle $\eN$ is $q$-ample, 
$q=\dim\mbb P(\eN^\vee)-\mbb P^{\nu+u}=\dim X-(u+1).$ Hence $Y=\Grs(\nu+u;\nu)$, the zero locus of a section of $\eN$, is $(u+1-\nu)^\pos$; this can be negative. 

On the other hand, Proposition~\ref{prop:x-b} implies that $Y$ is $u^\pos$. Moreover, $s_\ell$ above is far from being a regular section, so one can not use \ref{prop:q21} to estimate the amplitude of $Y_\ell$. But \ref{prop:x-b} is still applicable.\smallskip

\nit(iii) 
Zero loci of sections in globally generated vector bundles appear in recent work of Fulger-Lehmann (cf.~\cite{ful-leh}). They defined the pliant cone of a projective variety $X$, which is generated by pre-images of Schubert subvarieties of (various) Grassmannians $\Grs$, by morphisms $X\to\Grs$. The pliant cone is a full-dimensional sub-cone of the nef cone of $X$, whose generators are easier to understand. 

Pre-images of special Schubert varieties belong to the pliant cone. Our discussion shows they enjoy remarkable positivity properties, in particular are G3 in the ambient space.
\end{m-remark}


\subsection{Sources of $G_m$-actions}\label{ssct:fixed} 

Let $X$ be a smooth projective variety with a faithful action 
$$
\lda:G_m\times X\to X
$$ 
of the multiplicative group $G_m=\kk^\times$. This determines the so-called Bialynicki-Birula---BB for short---decomposition of $X$ (cf. \cite{bb}): 
\begin{itemize}[leftmargin=5ex]
	\item[$\bullet$] 
		The fixed locus $X^\lda$ of the action is a disjoint union 
		$\underset{s\in S_\BB}{\bigsqcup}\kern-1exY_s$ of smooth subvarieties. 
		For $s\in S_\BB$, 
		$Y_s^+:=\{x\in X\mid\underset{t\to 0}{\lim}\,\lda(t, x)\in Y_s\}$ 
		is locally closed in $X$ (a BB-cell) and it holds:
		$\;X=\underset{s\in S_\BB}{\bigsqcup}\kern-1exY_s^+.$
	\item[$\bullet$] 
		The \emph{source} $Y:=Y_{\rm source}$ and 
		the \emph{sink} $Y_{\text{\rm sink}}$ of the action are uniquely characterized 
		by the conditions: $Y^+=Y_{\rm source}^+\subset X$ is open 
		and $Y_{\rm sink}^+=Y_{\rm sink}$. 
\end{itemize}

A linearization of the action in a sufficiently ample line bundle yields a $G_m$-equivariant embedding $X\subset\mbb P^N_\kk$. There are homogeneous coordinates $\bz_{0}\in\kk^{N_0+1},\ldots,\bz_{r}\in\kk^{N_r+1}$ such that the $G_m$-action on $\mbb P^N_\kk$ is: 
\begin{equation}\label{eq:c*}
\lda\big(t,[\bz_{0},\bz_{1},\ldots,\bz_{r}]\big)
=[\bz_{0},t^{m_1}\bz_{1},\ldots,t^{m_r}\bz_{r}],\quad\text{with}\;0<m_1<\ldots<m_r. 
\end{equation}
The source and sink of $\mbb P^N, X$ are respectively: 
\begin{equation}\label{eq:YP}
\begin{array}{lcll}
\mbb P^N_{\rm source}=\{[\bz_{0},0,\ldots,0]\}, 
&&
\mbb P^N_{\rm sink}=\{[0,\ldots,0,\bz_{r}]\},
\\[1ex] 
Y=Y_{\rm source}=X\cap \mbb P^N_{\rm source},
&& 
Y_{\rm sink}=X\cap \mbb P^N_{\rm sink},
\\[1ex]
Y^+=X\cap (\mbb P^N_{\rm source})^+,
&&
(\mbb P^N_{\rm source})^+
=\{[\bz]=[\bz_{0},\bz_{1},\ldots,\bz_{r}]\mid\bz_{0}\neq 0\}.
\end{array}
\end{equation}

Let $m$ be the lowest common multiple of $\{m_\rho\}_{\rho=1,\dots,r}$ and $l_\rho:=m/m_\rho$. Denote $\bz_{\rho}^{l_\rho}:=(z_{\rho 0}^{l_\rho},\dots,z_{\rho N_\rho}^{l_\rho})$ and $\eI\subset\eO_{\mbb P^N}$ the sheaf of ideals generated by $\bz_1^{l_1},\dots,\bz_r^{l_r}$. The rational map
\begin{equation}\label{eq:phi}
\phi:\mbb P^N\dashto\mbb P^{N'},\quad
[\bz_0,\bz_1,\dots,\bz_r]\mt[\bz_1^{l_1},\dots,\bz_r^{l_r}],
\end{equation}
is $G_m$-invariant and its indeterminacy locus is the subscheme determined by $\eI$. Then $\eJ:=\eI\otimes\eO_{X}$ defines the subscheme $Y_\eJ\subset X$ whose reduction is $(Y,\eO_Y)$. We have the diagram: 
\begin{equation}\label{eq:blup}
\xymatrix@C=4em@R=2em{
\tld X:=\Bl_{Y_\eJ}(X)\ar[r]^-{\tld\iota}\ar[d]_-{b}
\ar@/^4ex/@<+.5ex>[rr]|{\;\phi_X\;}
&\Bl_\eI(\mbb P^N)\ar[r]^-{\phi}\ar[d]^-{B}
&\mbb P^{N'}\ar@{=}[d]
\\ 
X\ar[r]^-\iota&\mbb P^N\ar@{-->}[r]^-\phi&\mbb P^{N'}
}
\end{equation}

\begin{m-lemma}\label{lm:B}
The diagram \eqref{eq:blup} has the following properties: 
\begin{enumerate}[leftmargin=5ex]
\item[\rm(i)] 
The exceptional divisor of $B$ is $\phi$-relatively ample, hence the exceptional divisor of $b$ is $\phi_X$-relatively ample. 
\item[\rm(ii)] 
The morphism 
$\phi:\Bl_{\eI}(\mbb P^{N})\to\mbb P^{N'}$ is $G_m$-invariant and 
\begin{equation}\label{eq:estim}
\dim\phi_X(\tld X)=\dim\phi_X(X\sm Y^+)\les\dim (X\sm Y^+).
\end{equation}
\end{enumerate}
\end{m-lemma}

\begin{m-proof}
(i) The subscheme determined by $\eI$ is the vanishing locus of a section in a direct sum of ample line bundles over $\mbb P^N$, so we recover the situation in Proposition~\ref{prop:p}. 

\nit(ii) The $G_m$-invariance of $\phi$, thus of $\phi_X$, follows from \ref{eq:c*}. It holds: 
$$
\dim\phi_X(\Bl_\eJ(X))=\dim\ovl{\phi_X(X\sm Y)}
\;\;\text{and}\;\;
\phi_X(X\sm Y)=\phi_X(X\sm Y^+)\cup\phi_X(Y^+\sm Y).
$$ 
For $[\bz_0,\bz']\in Y^+\sm Y$ and $t\in G_m$, the $G_m$-invariance of $\phi_X$ yields:  
$$
\phi_X\big([\bz_0,\bz']\big)=\phi_X\big(t\times[\bz_0,\bz']\big)
=\phi_X\big(\lim_{t\to\infty}t\times[\bz_0,\bz']\big).
$$
But $\disp\lim_{t\to\infty}t\times[\bz_0,\bz']=[0,\bz'']\in X\sm Y^+$, which implies $\phi_X(Y^+\sm Y)\subset\phi_X(X\sm Y^+)$. 
\end{m-proof}

Now we can estimate the ampleness of the source $Y$. 

\begin{m-theorem}\label{thm:y-sink}
Let $X$ be a non-singular $G_m$-variety with source $Y$ and $$p:=\codim(X\sm Y^+)-1.$$ 
Then the following statements hold: 
\begin{enumerate}[leftmargin=5ex]
\item[\rm(i)] 
The thickening $Y_\eJ$ of $Y$ in \eqref{eq:blup} is a $(\dim Y-p)$-ample subscheme of $X$; in particular, $Y$ is a $p^\apos$ subvariety. 
\item[\rm(ii)]
If $G_m$ acts on the normal bundle $\eN_{Y/X}$ by scalar multiplication, then $Y\subset X$ is $p^\pos$.
\end{enumerate}  
\end{m-theorem}

\begin{m-proof}
(i) We apply the Proposition~\ref{prop:x-b}: $Y_\eJ$ is a $q$-ample subscheme, with
$$
q=1+\dim\phi(\tld X)-\codim_X(Y)\srel{\eqref{eq:estim}}{\les}
1+\dim(X\sm Y^+)-\codim_X(Y).
$$
\nit(ii) In this case we have $Y^+\cong\unbar{\sf N}:= \Spec\big(\Sym^\bullet\eN_{Y/X}^\vee\big)$, cf. \cite[Remark pp. 491]{bb}. Thus $\unbar{\sf N}\subset X$ is open and $G_m$ acts, fibrewise over $Y$, by scalar multiplication. 

The inclusions 
$\unbar{\sf N}\subset\unbar{\sf N}_{\mbb P^N_{\rm source}/\mbb P^N}
=
\{\,[\bz_{N_0},\bz']\mid\bz_{N_0}\neq0\,\}\subset\mbb P^N$ 
are $G_m$-equivariant. But the diagonal multiplication on the coordinates $\bz'$ exists globally on $\mbb P^N$, consequently $X\subset\mbb P^N$ is invariant for the $G_m$-action by scalar multiplication on $\bz'=(z_{N_1},\ldots,z_{N_r})$. Hence the exponents $l_\rho$ in \eqref{eq:phi} are all equal to one, so $\eJ=\eI_Y\subset\eO_X$. 
\end{m-proof}

\begin{m-remark}\label{rmk:cd}
Note that $X\sm Y^+$ is closed in $X\sm Y$. Thus, by using Lemma~\ref{lm:cd-apos}, we deduce: $\;\cd(X\sm Y)=\dim(X\sm Y^+)$. This simple answer contrasts the elaborate techniques needed to estimate the cohomological dimension of the complement of a subvariety (cf. \cite{ogus,falt-homog,lyub-licht}).
\end{m-remark}

\begin{m-example}\label{expl:o-gr} 
\nit(i) 
Let $W\cong\kk^{w+1}$, $w+1$ even, be a vector space endowed with a non-degenerate, symmetric bilinear form $\beta$. Consider $X:=\oGrs(u+1;W)$, the orthogonal Grassmannian of $(u+1)$-dimensional isotropic subspaces of $W$; in particular, $w+1\ges2(u+1)$. We choose coordinates on $W$ such that 
$$
\beta=\left[\begin{array}{cc}0&\bone_{(w+1)/2}\\ \bone_{(w+1)/2}&0\end{array}\right],
\quad\text{($\bone$ stands for the identity matrix)},
$$ 
and decompose $W=\kk^{(w+1)/2}\oplus\kk^{(w+1)/2}$ into the sum of two Lagrangian subspaces. We consider   
$\;\lda:G_m\to\SO_{(w+1)/2},\quad\lda(t)=\diag\big[t^{-1},\bone_{(w-1)/2},t,\bone_{(w-1)/2}\big].$

The source is $Y=\{U\mid s:=(1,0,\ldots,0)\in U\}\cong\oGrs(u;w-1)$ and, for $U\in Y$, holds: 
$$
\begin{array}{rl}
\eT_{X,U}
&\cong\Hom(U,U^\perp/U)\oplus\Hom^{\rm anti-symm}(U,U^\vee),
\qquad\text{(Note that $h(s)\in\lran{s}^\perp$.)}
\\[1ex] 
\eT_{Y,U}
&
\cong\Hom(U/\lran{s},U^\perp/U)\oplus
\Hom^{\rm anti-symm}\big(U/\lran{s},(U/\lran{s})^\vee\big), 
\\[1ex] 
\eN_{Y/X, U}
&=\Hom(\lran{s},\lran{s}^\perp/U). 
\end{array}
$$
Hence $\lda$ acts with weight $t$ on $\eN_{Y/X}$. The complement of the open BB-cell is 
$$
X\sm Y^+=\big\{U\in X\mid s\not\in\underset{t\to 0}{\lim}\lda(t)U\big\}
=
\{U\mid U\subset W':=\kk^{(w-1)/2}\oplus\kk^{(w+1)/2}\}.
$$ 
Note that $\beta_{\rst_{W'}}$ has a $1$-dimensional kernel $\lran{s'}$: if $w=2u+1$, then $s'\in U$ for all $U\in X\sm Y^+$; for greater $w$, this is not the case. We deduce: 
$$
\codim(X\sm Y^+)=\bigg\{ 
\begin{array}{cl}
u&\text{if } w=2u+1;\\ u+1&\text{if }w\ges2u+3,
\end{array}
\;\Rightarrow\;
\text{$Y\subset X$ is: } 
\bigg\{
\begin{array}{cl}
(u-1)^\pos&\text{if }w=2u+1;\\ u^\pos&\text{if }w\ges2u+3.
\end{array}
$$ 

\nit(ii) 
With the previous notation, let $\omega$ be a skew-symmetric bilinear form on $W$ and $X:=\spGrs(u+1;W)$ be the symplectic Grassmannian of $(u+1)$-dimensional isotropic subspaces of $W$. Take a Lagrangian decomposition $W=\kk^{(w+1)/2}\oplus\kk^{(w+1)/2}$, such that 
$$
\omega=\left[\begin{array}{cc}0&\bone_{(w+1)/2}\\ -\bone_{(w+1)/2}&0\end{array}\right].
$$
The action induced by  $\;\lda:G_m\to\Sp_{(w+1)/2},$ $\lda(t)=\diag\big[t^{-1},\bone_{(w-1)/2},t,\bone_{(w-1)/2}\big]$ has the source 
$Y=\{U\mid s:=(1,0,\ldots,0)\in U\}\cong\spGrs(u;w-1)$. 

A similar computation yields  $\eN_{Y/X,U}\cong\Hom(\lran{s},W/U).$ In this case, $G_m$ doesn't act on $\eN_{Y/X}$ by scalar multiplication: it has weight $t^2$ on $\Hom(\lran{s},W/\lran{s}^\perp)$ and weight $t$ on its complement ($\perp$ stands for the $\omega$-orthogonal). The complement of the open BB-cell is $(1+u)$-codimensional in $X$, as before. We conclude that, in this case, $Y$ is only $u^\apos$; more precisely, there is a non-reduced scheme with support $Y$ which is $u^\pos$. 
\end{m-example}


\end{document}